\documentclass[11pt]{amsart}
\usepackage{mathrsfs, amsmath, amssymb, array,graphics, booktabs,diagbox,listings,color,textcomp,picinpar,wrapfig,cutwin,longtable,slashbox,amsthm,tikz,threeparttable,makecell,cases}

\usepackage[backref]{hyperref}
\usepackage{multicol}
\usepackage{multirow}
\usepackage[page,title,titletoc,header]{appendix}
\usepackage[font=small,labelfont=bf]{caption}
\DeclareCaptionFont{tiny}{\tiny}
\DeclareCaptionFont{normalsize}{\normalsize}

\usepackage{geometry}
\definecolor{listinggray}{gray}{0.9}

\definecolor{lbcolor}{rgb}{0.9,0.9,0.9}

\usepackage{algorithm}
\usepackage[noend]{algpseudocode}

\newtheorem{theorem}{Theorem}[section]

\newtheorem{corollary}[theorem]{Corollary}
\newtheorem{proposition}[theorem]{Proposition}

\theoremstyle{definition}
\newtheorem{definition}[theorem]{Definition}

\theoremstyle{remark}
\newtheorem{remark}[theorem]{Remark}

\numberwithin{equation}{section}

\newcommand{\HH}{\mathbb{H}}

\begin{document}

\title{Compact hyperbolic Coxeter five-dimensional polytopes with nine facets}

\author{Jiming Ma}
\address{School of Mathematical Sciences \\Fudan University\\Shanghai 200433, China} \email{majiming@fudan.edu.cn}

\author{Fangting Zheng}
\address{Department of Mathematical Sciences \\ Xi'an Jiaotong Liverpool University\\ Suzhou 200433,	China  }
\email{Fangting.Zheng@xjtlu.edu.cn}

\date{Nov. 22, 2022}
\thanks{Jiming Ma was partially supported by  NSFC  12171092. Fangting Zheng was supported by NSFC 12101504 and XJTLU Research Development Fund RDF-19-01-29}

\begin{abstract}
	In this paper, we obtain a complete classification of compact hyperbolic Coxeter five-dimensional polytopes with nine facets.
\end{abstract}

\maketitle

\section{Introduction}

A Coxeter polytope in the spherical, hyperbolic or Euclidean space is a polytope whose dihedral angles are all integer submultiples of $\pi$. Let $\mathbb{X}^d$ be $\mathbb{E}^d$, $\mathbb{S}^d$, or $\mathbb{H}^d$. If $\Gamma\subset Isom(\mathbb{X}^d)$ is a finitely generated discrete reflection group, then its fundamental domain is a Coxeter polytope in $\mathbb{X}^d$. On the other hand, if $\Gamma=\Gamma(P)$ is generated by reflections in the bounding hyperplanes of a Coxeter polytope $P\subset\mathbb{X}^d$, then $\Gamma$ is a discrete group of isometries of $\mathbb{X}^d$ and $P$ is its fundamental domain.

There is an extensive body of literature in this field. In early work, \cite{Coxeter:1934} has proved that any spherical Coxeter polytope 
 is a simplex and any Euclidean Coxeter polytope is either a simplex or a direct product of simplices. See, for example, \cite{Coxeter:1934,Bourbaki: 1968}  for full lists of spherical and Euclidean Coxeter polytopes.

However, for hyperbolic Coxeter polytopes, the classification remains an active research topic. It was proved by Vinberg \cite{Vinberg:1985} that no compact hyperbolic Coxeter polytope exists in dimensions $d\geq 30$, and non-compact hyperbolic Coxeter polytope of finite volume does not exist in dimensions $d \geq 996$ \cite{Prokhorov:1987}. These bounds, however, may not be sharp. Examples of compact polytopes are known up to dimension $8$ \cite{Bugaenko:1984,Bugaenko:1992}; non-compact polytopes of finite volume are known up to dimension $21$ \cite{Vinberg:1972,VK:1978,Borcherds: 1998}. As for the classification, complete results are only available in dimensions less than or equal to three. Poincare finished the classification of 2-dimensional hyperbolic polytopes in \cite{Poincare: 1882}. That result was important to the work of Klein and Poincare on discrete groups of isometries of the hyperbolic plane. In 1970, Andreev proved an analogue for the $3$-dimensional hyperbolic convex finite volume polytopes \cite{Andreev1: 1970,Andreev2: 1970}. This theorem played a fundamental role in Thurston's work on the  geometrization of 3-dimensional Haken manifolds.

In higher dimensions, although a complete classification is not available, interesting examples have been displayed in \cite{Makarov: 1965,Makarov: 1966,Vinberg: 1967,Makarov: 1968,Vinberg: 1969,Rusmanov: 1989,ImH: 1990,Allcock: 2006}. In addition, enumerations are reported for the cases in which the differences between the numbers of facets $m$ and the  dimensions $d$ of polytopes are fixed to some small number. When $m-d=1$, Lann\'{e}r classified all compact hyperbolic Coxeter simplices \cite{Lanner: 1950}. The enumeration of non-compact hyperbolic simplices with finite volume has been reported by several authors, see e.g. \cite{Bourbaki: 1968,Vinberg: 1967,Koszul:1967}. For $m-d=2$, Kaplinskaja described all hyperbolic Coxeter simplicial prisms \cite{Kaplinskaya: 1974}. Esselmann \cite{Esselmann: 1996} later enumerated other possibilities in this family, which are named \emph{Esselmann polytopes}. 
Tumarkin \cite{Tumarkin: n2} classified all non-compact hyperbolic Coxeter $d$-dimensional polytopes with
$n + 2$ facets. In the case of $m-d=3$, Esselman proved in 1994 that compact hyperbolic Coxeter $d$-polytopes with $d+3$ facets  only exist when $d\leq 8$ \cite{Esselmann}. By expanding the techniques derived by Esselmann in \cite{Esselmann} and \cite{Esselmann: 1996}, Tumarkin completed the classification of compact hyperbolic Coxeter $d$-polytopes with $d+3$ facets \cite{Tumarkin: n3}. In the non-compact case, Tumarkin proved in \cite{Tumarkin: n3fv,Tumarkin: n3fvs} that such polytopes do not exist in dimensions greater than or equal to $17$. Furthermore, Roberts provided a list for the family with exactly one non-simple vertex \cite{Roberts:15}. In the case of $m-d=4$,
Flikson-Tumarkin showed in \cite{FT:08} that compact hyperbolic Coxeter $d$-polytope with $d + 4$ facets does not exist when $d$ is greater than or equal to $8$. This bound is sharp because of the example constructed by Bugeanko \cite{Bugaenko:1984}. In addition, Flikson-Tumarkin showed that Bugeanko's example is the only $7$-dimensional polytope with $11$ facets. However, complete classifications for $d=4,5,6$ are not presented. Recently,  Burcroff \cite{Amanda:2022} and Ma-Zheng \cite{MZ:2022} provided a complete list of $348$ compact Coxeter $4$-polytopes with $8$ polytopes independently. 

Besides, some scholars have also considered polytopes with small numbers of disjoint pairs
\cite{FT:08s,FT:09,FT:14} or of  certain combinatoric types, such as $d$-pyramid \cite{Tumarkin: n2,Tumarkin: n3fv} and $d$-cube \cite{Jacquemet2017,JT:2018}. An updating overview of the current knowledge for hyperbolic Coxeter polytopes is available on Anna Felikson's webpage \cite{Annahomepage}.

In this paper, the following is proved:

\begin{theorem}\label{thm:main}
	There are exactly $51$ compact hyperbolic Coxeter $5$-polytopes with $9$ facets. 
\end{theorem}

We remark that Burcroff has also carried out $50$ such polytopes independently \cite{Amanda:2022}. We have communicated with Burcroff our results of $51$ polytopes when her preprint appeared. Burcroff replied that she lost the case where the weights are all less than or equal to five. We now all agree that $51$ is the correct number. The correspondence between the notions of our and Burcroff’s polytopes is presented in Section Section \ref{section:vadilation}.

The paper \cite{JT:2018} is the main inspiration for our recent work on enumerating hyperbolic Coxeter polytopes. In comparison with \cite{JT:2018}, we use a more universal ``block-pasting" algorithm, which is first introduced in \cite{MZ:2018}, rather than the ``tracing back" algorithm. More geometric obstructions are adopted and programmed to considerably reduce the computational load. Our algorithm efficiently enumerates hyperbolic Coxeter polytopes over arbitrary combinatoric type rather than merely the $n$-cube. 

Last but not the least, our main motivation in studying the hyperbolic Coxeter polytopes is for  the construction of high-dimensional hyperbolic manifolds. However, this is not the theme here. Readers can refer to, for example,  \cite{KM: 2013}, for interesting hyperbolic manifolds built via special hyperbolic Coxeter polytopes.

 This  paper is organized as follows. We provide in Section $2$ some preliminaries about hyperbolic (compact Coxeter) polytopes. Then, in Section $3$ we recall the 2-phases procedure and some terminologies introduced by Jacquemet and Tschantz \cite{Jacquemet2017,JT:2018} for numerating hyperbolic Coxeter $n$-cubes. The $322$ combinatorial types of simple $5$-polytopes with $9$ facets are reported in Section $4$. The enumerations of all the ``SEILper"-potential matrices are in Section $5$. The sectional restrictions are also introduced there to further restrict the number of SEILper matrices. Next, signature  obstructions are applied for Gram matrices of actual hyperbolic Coxeter polytopes in Section $6$. Validations and the complete lists of the resulting Coxeter diagrams and hyperbolic lengths of Theorem \ref{thm:main} are shown in Section $7$. 

\textbf{Acknowledgment}

We would like to thank Amanda Burcroff for communicating with us about her result after we posted our preprint \cite{MZ:2022}. The computations is pretty delicate and complex, and the list now is much more  convincing due to the  mutual check. We are also grateful to Nikolay Bogachev for his interest and discussion about the results. The computations throughout this paper are performed on a cluster of server of PARATERA, engrid12, line priv$\_$para (CPU:Intel(R) Xeon(R) Gold 5218 16 Core v5@2.3GHz).

\section{preliminary} \label{section:cchp}
In this section, we recall some essential facts about compact Coxeter hyperbolic polytopes, including Gram matrices, Coxeter diagrams, characterization theorems, etc. Readers can refer to, for example, \cite{Vinberg:1993} for more details.  

\subsection{Hyperbolic space, hyperplane and convex polytope}
We first describe a hyperboloid model of the $d$-dimensional hyperbolic space $\mathbb{H}^d$. Let $\mathbb{E}^{d,1}$ be a $d+1$-dimensional Euclidean vector space equipped with a Lorentzian scalar product $\langle\cdot,\cdot \rangle$ of signature $(d,1)$. We denote by $C_+$ and $C_-$ the connected components of the open cone $$C=\{x=(x_1, ...,x_d,x_{d+1})\in \mathbb{E}^{d,1}:\langle x,x\rangle<0\}$$ 
with $x_{d+1}>0$ and $x_{d+1}<0$, respectively. Let $R_{+}$ be the group of positive numbers acting on $\mathbb{E}^{d,1}$ by homothety. The hyperbolic space $\mathbb{H}^d$ can be identified with the quotient set $C_+/R_+$, which is a subset of $P\mathbb{S}^d=(\mathbb{E}^{d,1}\backslash \{0\})/R_+.$ The natural projection is denoted by $$\pi:(\mathbb{E}^{d,1}\backslash \{0\})\rightarrow P\mathbb{S}^d.$$ 

\noindent We denote $\overline{\mathbb{H}^d}$ as the completion of $\mathbb{H}^d$ in $P\mathbb{S}^d$. The points of the boundary $\partial \mathbb{H}^d= \overline{\mathbb{H}^d}\backslash \mathbb{H}^d$ are called the \emph{ideal points}. 
The affine subspaces of $\mathbb{H}^d$ of dimension $d-1$ are \emph{hyperplanes}. In particular, every hyperplane of $\mathbb{H}^d$ can be represented as $$H_e=\{\pi(x):x\in C_+,\langle x,e\rangle=0\},$$ where $e$ is a vector with $\langle e,e \rangle=1$. The half-spaces separated by $H_e$ are denoted by $H_e^+$ and $H_e^-$, where
\begin{equation}
H_e^-=\{\pi(x):x\in C_+,\langle x,e\rangle\leq 0\}. \label{1}
\end{equation} 
  
The \emph{mutual disposition of hyperplanes} $H_{e}$ and $H_{f}$ can be described in terms of the corresponding two vectors $e$ and $f$ as follows: 
\begin{itemize}
	\item The hyperplanes $H_{e}$ and $H_{f}$ intersect if $\vert\langle e,f \rangle\vert<1$. The value of the dihedral angle of $H_{e}^-\cap H_{f}^-$, denoted by $\angle H_e H_f$, can be obtained via the formula $$\cos \angle H_e H_f=-\langle e,f\rangle;$$
	\item The hyperplanes $H_{e}$ and $H_{f}$ are ultraparallel if  $\vert\langle e,f\rangle\vert=1$;
    \item The hyperplanes $H_{e}$ and $H_{f}$ diverge if $\vert\langle e,f\rangle\vert>1$. The distance $\rho(H_e,H_f)$ between $H_e$ and $H_f$, when $H_e^+\subset H_f^-$ and $H_f^+\subset H_e^-$, is determined by  $$\cosh \rho(H_e,H_f)=-\langle e,f \rangle.$$
\end{itemize}

 
We say a hyperplane $H_e$ \emph{supports} a closed bounded convex set $S$ if $H_e\cap S \ne \empty 0$ and $S$
lies in one of the two closed half-spaces bounded by $H_e$. If a hyperplane $H_e$ supports $S$, then $H_e \cap S$ is called a \emph{face} of $S$. 

\begin{definition}
A $d$-dimensional convex hyperbolic polytope is a subset of the form 

\begin{equation}
P=\overline{\mathop{\cap}\limits_{i\in \mathcal{I}} H_{i}^-}\in\overline{\mathbb{H}^d}, \label{2}
\end{equation}
where $H_i^-$ is the negative half-space bounded by the hyperplane $H_i$ in $\mathbb{H}^d$ and the line ``--- "  above the intersection  means taking the completion in $\overline{\mathbb{H}^d}$, 
under the following assumptions:
\begin{itemize}
	\item  $P$ contains a non-empty open subset of $\mathbb{H}^d$ and is of finite volume;
	\item  Every bounded subset $S$ of $P$ intersects only finitely many $H_{i}$.
\end{itemize}
\end{definition}

A convex polytope of the form (\ref{2}) is called \emph{acute-angled} if for distinct $i,j$, either $\angle H_iH_j\leq \frac{\pi}{2}$ or $H_i^+\cap H_j^+=\emptyset$. It is obvious that Coxeter polytopes are acute-angled. We denote  $e_i$ as the corresponding unit vector of $H_i$, namely $e_i$ is orthogonal to $H_i$ and point away from $P$. The polytope $P$ has the following form in the hyperboloid model.
 $$P=\pi(K)\cap \overline{\mathbb{H}^d},$$ where $K=K(P)$ is the convex polyhedral cone in $\mathbb{E}^{d,1}$ given by $$K=\{x\in \mathbb{E}^{d,1}: \langle x,e_i\rangle \leq 0~ \text{for all}~ i \}.$$ 

In the sequel, a $d$-dimensional convex polytope $P$ is called a \emph{$d$-polytope}. A $j$-dimensional face is named a $j$-face of $P$. In particular, a $(d-1)$-face is called a \emph{facet} of $P$.  We assume that each of the hyperplane $H_i$ intersects with $P$ on its facet. In other words, the hyperplane $H_i$ is uniquely determined by $P$ and is called a \emph{bounding hyperplane} of the polytope $P$. A hyperbolic polytope $P$ is called \emph{compact} if all of its $0$-faces, i.e., vertices, are in $\mathbb{H}^d$. It is called of  \emph{finite volume} if some vertices of $P$ lie in $\partial\mathbb{H}^d$.

\subsection{Gram matrices, Perron-Frobenius Theorem, and Coxeter diagrams}\label{hcp}

Most of the content in this subsection is well-known by  peers in this field. We present them for the convenience of the readers. 
In particular, Theorem \ref{thm:signature} and \ref{Vinberg:thm3.1} are extremely important throughout this paper. 

For a hyperbolic Coxeter $d$-polytope $P=\overline{\mathop{\cap}\limits_{i\in \mathcal{I}} H_{i}^-}$, we define the Gram matrix of polytope $P$ to be the Gram matrix $(\langle e_i,e_j\rangle)$ of the system of vectors $\{e_i\in \mathbb{E}^{d,1}\vert i\in\mathcal{I}\}$ that determine $H_i^-$s. 
The Gram matrix of $P$ is the $m\times m$ symmetric matrix $G(P)=(g_{ij})_{1\leq i,j\leq m}$  defined as follows:

\begin{center}
	$g_{ij}=
	\left\{
	\begin{array}{ccl}1  \hspace{0.7cm} &\mbox~~~~~{{\rm if}}& j=i, \\
	-\cos\frac{\pi}{k_{ij}}&\mbox ~~~~~{{\rm if}}&  H_i~ \text{and}~ H_j~ \text{intersect} \text{ at~a~dihedral~ angle~}~ \frac{\pi}{k_{ij}},\\
	-\cosh \rho_{ij}&\mbox ~~~~~{{\rm if}}&  H_i~ \text{and}~ H_j~ \text{divergeat~at~a~distance}~\rho_{ij},\\
	-1  \hspace{0.7cm} &\mbox~~~~~{{\rm if}}& H_i~\text{and}~H_j~\text{are~ultraparallel}.
	\end{array}\right.
	$
\end{center}

 Other than the Gram matrix,  a Coxeter polytope $P$ can also be described by its \emph{Coxeter graph} $\Gamma=\Gamma(P)$. Every node $i$ in $\Gamma$ represents the bounding hyperplane $H_i$ of $P$. Two nodes $i_1$ and $i_2$ are joined by an edge with weight $2\leq k_{ij} \leq \infty$ if $H_i$ and $H_j$ intersect in $\mathbb{H}^n$ with angle $\frac{\pi}{k_{ij}}$. If the hyperplanes $H_i$ and $H_j$ have a common perpendicular of length $\rho_{ij}>0$ in $\mathbb{H}^n$, the nodes $i_1$ and $i_2$ are joined by a dotted edge, sometimes labelled $\cosh \rho_{ij}$. In the following, an edge of weight $2$ is omitted, and an edge of weight $3$ is written without its weight. The rank of $\Gamma$ is defined as the number of its nodes. In the compact case, $k_{ij}$ is not $\infty$, and we have $2\leq k_{ij}< \infty$. 
 
 A square matrix $M$ is said to be the direct sum of the matrices $M_1, M_2,\cdots,M_n$ if by some permutation of the rows and of columns, it can be brought to the form
 \begin{center}
 	$ \begin{pmatrix}
 	M_1&&&&0 \\
 	&M_2&&&\\
 	&& \cdot&&\\
 	&&&\cdot& \\
 	0&&&&M_n
 	\end{pmatrix}_. $
 	
 \end{center}

 \noindent A matrix $M$ that cannot be represented as a direct sum of two matrices is said to be \emph{indecomposible}\footnote{It is also referred to as ``irreducible" in some references.}. Every matrix can be represented uniquely as a direct sum of indecomposible matrices, which are called (indecomposible) components. 
 We say a polytope is \emph{indecomposible} if its Gram matrix $G(P)$ is indecomposible.
 
 \begin{figure}[h]
 	\scalebox{0.35}[0.35]{\includegraphics {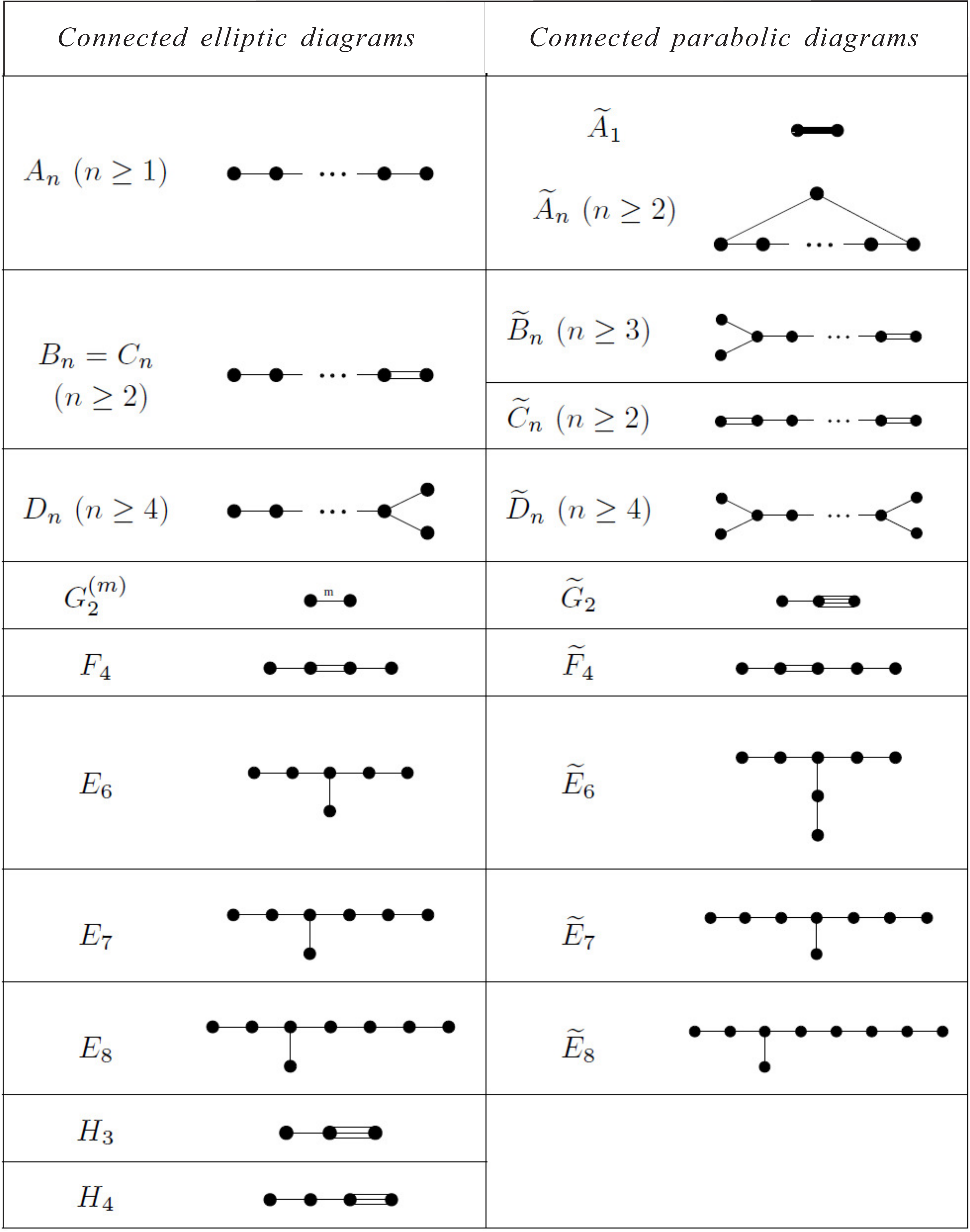}}
 	\caption{Connected elliptic (left) and connected parabolic (right) Coxeter diagrams.}
 	\label{figure:coxeter}
 \end{figure}
 
 In 1907, Perron found a remarkable property of the the eigenvalues and eigenvectors of matrices with positive entries. Frobenius later generalized it by investigating the spectral properties of indecomposible non-negative matrices.
 
 \begin{theorem}[Perron-Frobenius, \cite{G:1959}]\label{thm:PF}
 	An indecomposible matrix $A=(a_{ij})$ with non-positive entries always has a single positive eigenvalue $r$ of $A$. The corresponding eigenvector has positive coordinates. The moduli of all of the other eigenvalues do not exceed $r$.
 \end{theorem}

 It is obvious that Gram matrices $G(P)$ of an indecomposible Coxeter polytope is an indecomposible symmetric matrix with non-positive entries off the diagonal. Since the diagonal elements of $G(P)$ are all $1$s, $G(P)$ is either positive definite, semi-positive definite or indefinite. According to the Perron-Frobenius theorem, the defect of a connected semi-positive definite matrix $G(P)$ does not exceed $1$, and any proper submatrix of it is positive definite. For a Coxeter $n$-polytope $P$, its Coxeter diagram $\Gamma(P)$ is said to be \emph{elliptic} if $G(P)$ is positive definite; $\Gamma (P)$ is called \emph{parabolic} if any indecomposable component of $G(P)$ is degenerate and every subdiagram is elliptic. The elliptic and connected parabolic diagrams are exactly the Coxeter diagrams of spherical and Euclidean Coxeter simplices, respectively. They are classified  by Coxeter \cite{Coxeter:1934} as shown in Figure \ref{figure:coxeter}.
 
\begin{figure}[h]
 	\scalebox{0.35}[0.35]{\includegraphics {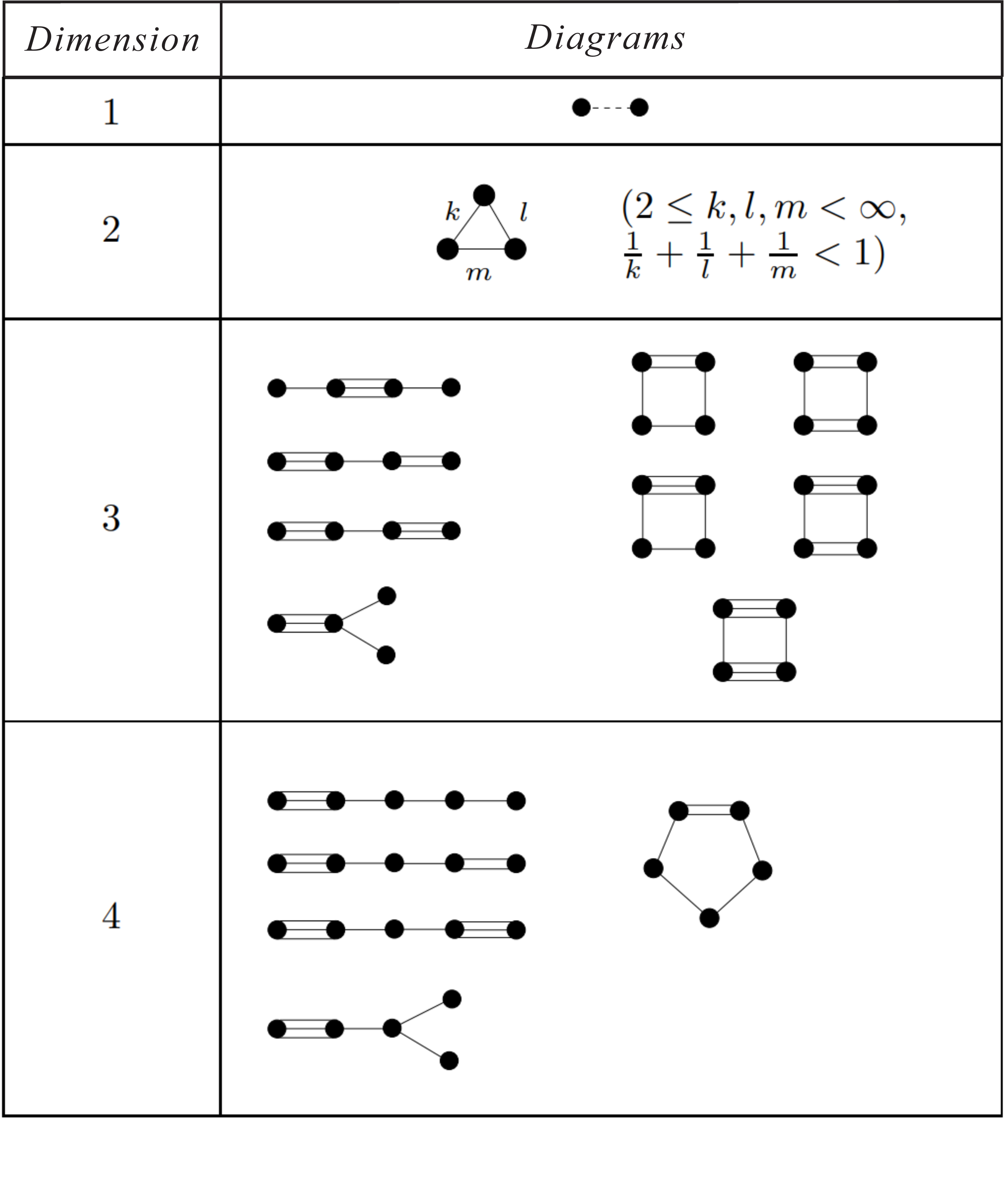}}
 	\caption{The Lann\'{e}r diagrams.}
 	\label{figure:lanner}
 \end{figure}

  A connected diagram $\Gamma$ is a \emph{Lann\'{e}r diagram} if $\Gamma$ is neither elliptic nor parabolic; any proper subdiagram of $\Gamma$ is elliptic. Those diagrams are irreducible Coxeter diagrams of compact hyperbolic Coxeter simplices. All such diagrams, reported by Lann\'{e}r \cite{Lanner: 1950}, are listed in Figure \ref{figure:lanner}.

Although the full list of hyperbolic Coxeter polytopes remains incomplete, some powerful algebraic restrictions are known \cite{Vinberg:1985}:
 
 \begin{theorem} \label{thm:signature}
 	(\cite{Vinberg:1985s}, Th. 2.1). Let $G=(g_{ij})$ be an indecomposable symmetric matrix
 	of signature $(d,1)$, where $g_{ii}=1$ and $g_{ij}\leq 0$ if $i\ne j$. Then there exists a unique (up to isometry of $\mathbb{H}^d$) convex hyperbolic polytope $P\subset\mathbb{H}^d$, whose Gram matrix coincides with $G$.
 \end{theorem}
 
 \begin{theorem} \label{Vinberg:thm3.1} 
 	(\cite{Vinberg:1985s}, Th. 3.1, Th. 3.2)
 	Let $P=\mathop{\cap}\limits_{i\in I} H_i^- \in \mathbb{H}^d$ be a compact acute-angled polytope and $G=G(P)$ be the Gram matrix. Denote $G_J$ the principal submatrix of G formed from the rows and columns whose indices belong to $J\subset I$. Then, 
 	\begin{enumerate}
 		\item The intersection $\mathop{\cap}\limits_{j\in J}H_j^-,J\subset I$, is a face $F$ of $P$ if and only if the matrix $G_J$ is positive definite
 		\item   For any $J\subset I$ the matrix,  $G_J$ is not parabolic.
 	\end{enumerate}
 \end{theorem}
 
  A convex polytope is said to be \emph{simple} if each of its faces of codimension $k$ is contained in exactly $k$ facets. 
  By Theorem \ref{Vinberg:thm3.1}, we have the following corollary:
 
 \begin{corollary}
 	Every compact acute-angled polytope is simple.
 \end{corollary}

\section{Potential hyperbolic Coxeter matrices}\label{section:potential}

In order to classify all of the compact hyperbolic Coxeter $5$-polytopes with $9$ facets, we firstly enumerate all Coxeter matrices of simple $5$-polytope with $9$ facets that satisfy spherical conditions around all of the vertices. These are named \emph{potential hyperbolic Coxeter matrices} in \cite{JT:2018}. Almost all of the terminology and theorems in this section are proposed by Jacquemet and Tschantz. We recall them here for reference, and readers can refer to \cite{JT:2018} for more details.

\subsection{Coxeter matrices}\label{subsection:coxeterMatrix}

The \emph{Coxeter matrix} of a hyperbolic Coxeter polytope $P$ is a symmetric matrix $M=(m_{ij})_{1\leq i,j\leq N}$ with entries in $\mathbb{N}\cup\{\infty\}$ such that
\[m_{ij}=\left\{\begin{array}{cl} 
1,&\text{if } j=i,\\ 
k_{ij}, &\text{if }H_i\text{ and }H_j\text{ intersect in }\HH^n\text{ with angle }\frac{\pi}{k_{ij}},\\
\infty, & \text{otherwise}.
\end{array}\right.\]
Note that, compared with Gram matrix, the Coxeter matrix does not involve the specific information of the distances of the disjoint pairs.

\begin{remark}
	In the subsequent discussions, we refer to \textit{the Coxeter matrix $M$ of a graph $\Gamma$} as the Coxeter matrix $M$ of the Coxeter polyhedron $P$ such that $\Gamma=\Gamma(P)$.\\
\end{remark}

\subsection{Partial matrices}\label{subsection:partial}

\begin{definition}
	Let $\Omega=\{n\in \mathbb{Z}\,|\,n\geq 2\}\cup\{\infty\}$ and let $\bigstar$ be a symbol representing an undetermined real value. A \textit{partial matrix of size $m\geq 1$} is a symmetric $m\times m$ matrix $M$ whose diagonal entries are $1$, and whose non-diagonal entries belong to $\Omega\cup\{\bigstar\}$.
\end{definition}

\begin{definition}
	Let $M$ be an arbitrary $m\times m$ matrix, and $s=(s_1,s_2,\cdots,s_k)$, $1\leq s_1<s_2<\cdots<s_k\leq m$. Let $M^{s}$ be the $k\times k$ submatrix of $M$ with $(i,j)$-entry $m_{s_i,s_j}$. 
\end{definition}

\begin{definition}
	We say that a partial matrix $M=(p_{ij})_{1\leq i,j,\leq m}$ is a \emph{potential matrix} for a given polytope $P$ if 
	
	$\bullet$ There are no entries with the value $\bigstar$;
	
	$\bullet$ There are entries $\infty$ in positions of $M$ that correspond to disjoint pair;
	
	$\bullet$ For every sequence $s$ of indices of facets meeting at a vertex $v$ of $P$, the matrix, obtained by replacing value $n$ with $\cos\frac{\pi}{n}$ of submatrix $M^s$, is elliptic.
	
\end{definition}

For brevity, we use a \emph{potential vector} $$C=(p_{12},p_{13},\cdots p_{1m},p_{23},p_{24},\cdots,p_{2m},\cdots p_{ij},\cdots p_{m-1,m}),~ p_{ij}\ne \infty$$ to denote the potential matrix, where $1\leq i<j\leq m$ and non-infinity entries are placed by the subscripts lexicographically. The potential matrix  and potential vector $C$  can be constructed one from each other easily. In general, an arbitrary Coxeter matrix corresponds to a Coxeter vector following the same manner. We mainly use the language of \emph{vectors} to explain the methodology and  report the enumeration results. It is worthy to remark that, for a given Coxeter diagram, the corresponding (potential / Gram) matrix and vector are not unique in the sense that they are determined under a given labelling system of the facets and may vary when the system changed. In Section \ref{chapter:algorithm}, we apply a permutation group to the nodes of the diagram and remove the duplicates to obtain all of the distinct desired vectors

For each rank $r\geq 2$, there are infinitely many finite Coxeter groups, because of the infinite 1-parameter family of all dihedral groups, whose graphs consist of two nodes joined by an edge of weight $k\geq 2$. However, a simple but useful truncation can be utilized:

\begin{proposition}
	There are finitely many finite Coxeter groups of rank $r$ with Coxeter matrix entries at most seven.
\end{proposition}

It thus suffices to enumerate potential matrices with entries at most seven, and the other candidates can be obtained from substituting integers greater than seven with the value seven. In other words, we now have more variables, that are restricted to be integers larger than or equal to seven, besides length unknowns.  In the following, we always use the terms ``Coxeter matrix" or ``potential matrix" to mean the one with integer entries less than or equal to seven unless otherwise mentioned.  

In \cite{JT:2018}, the problem of finding certain hyperbolic Coxeter polytopes is solved in two phases. In the first step, potential matrices for a particular hyperbolic Coxeter polytope are found; the ``Euclidean-square obstruction" is used to reduce the number. Secondly, relevant algebraic conditions are solved for the admissible distances between non-adjacent facets.  
In our setting, additional universal necessary conditions, except for the vertex spherical restriction and Euclidean square obstruction, are adopted and programmed to reduce the number of the potential matrices.

\section{ Simple 5-polytopes with 9 facets}\label{section:4d8f}
K. Fukuda, H. Miyata  and S. Moriyama completed the enumeration of small realizable oriented
matroids and, as a corollary, provided a list of $322$ $5$-dimensional simplicial polytopes with $9$ vertices in  \cite{FMM:13}, denoted by $P_1$, $P_2$, $\cdots$, $P_k$, $\cdots$, $P_{322}$, which is dual to the simple $5$-dimensional polytopes with $9$ facets. Each line in data file provided in \cite{FMM:13} is for one combinatorial type. The data on line $322$ is as follows:

{\color{blue}
\noindent [9,8,7,6,5] [9,8,7,6,4] [9,8,7,5,4] [9,8,6,5,4] [9,7,6,5,4] [8,7,6,5,3] [8,7,6,4,3] [8,7,5,4,3] [8,6,5,4,3] 
[7,6,5,4,2] [7,6,5,3,2] [7,6,4,3,2] [7,5,4,3,2] [6,5,4,3,1] [6,5,4,2,1] [6,5,3,2,1] [6,4,3,2,1] [5,4,3,2,1] 

  }

\noindent where the number $1$, $2$, $\cdots$, $9$ denote the nine facets and each square bracket corresponds to one vertex that is incident to the enclosed five facets. For example, there are $18$ vertices of the above polytope $P_{322}$.

From the original information, we can search out the following \textbf{\emph{data}}  for each polytope:

\begin{enumerate}
	\item The permutation subgroup $g_i$  of $S_9$ that is isomorphic to the symmetry group of $P_k$;
	\item The set $d_k$ of pairs of the disjoint facets;
	\item The set $l_4$ / $l_5$ of sets of four / five facets of which the intersection is of the combinatorial type of a tetrahedron / $4$-simplex.
	\item The set $l_5\_basis$ / $l_4\_basis$  of sets of five / four facets of which bound a 4-simplex / 3-simplex facet. The label of the bounded 4-somplex / 3-simplex facet is recorded as well. Note that $l_5\_basis$ / $l_4\_basis$  is a subset of $l_5$ / $l_4$. The set $l_j\_basis$ can be non-empty only for a $j$-dimensional polytope.
	\item The set $i_2$ of sets of facets of which the intersection is of the combinatorial type of a $2$-cube.
	\item The set $s_3~ /~ s_4 /~s_5 $ of sets of three / four / five facets of which the intersection is not a face /an edge / a vertex of $P_k$, and no disjoint pairs are included.
	\item The set $e_3~ /~ e_4 /~e_5$ of sets of three / four / five facets of which no disjoint pairs are included.
	\item The set $se_6~ /~ se_7$ of sets of six/ seven facets of which no disjoint pairs are included.
	
\end{enumerate}  

For example, for the polytope $P_{322}$, the above sets are as shown in Table \ref{table:p322}.

\begin{table}[h]
	{\footnotesize
		\begin{tabular}{|c|c|l|}
			\Xcline{1-3}{1.2pt}
			\multicolumn{3}{|c|}{\textbf{$P_{322}$ }}\\
			\hline
			\multirow{4}{*}{Vert}& \multirow{4}{*}{18}&  $\{	\{1, 2, 3, 4, 5\}, \{1, 2, 3, 4, 6\}, \{1, 2, 3, 5, 6\}, \{1, 2, 4, 5, 6\}, \{1, 3, 4, 5, 6\},\{2, 3, 4, 5, 7\}, \{2, 3, 4, 6, 7\}, $\\
			&&$	\{2, 3, 5, 6, 7\}, \{2, 4, 5, 6, 7\}, \{3, 4, 5, 6, 8\}, \{3, 4, 5, 7, 8\}, \{3, 4, 6, 7, 8\}, \{3, 5, 6, 7, 8\}, \{4, 5, 6, 7, 9\},$\\
			&&$  \{4, 5, 6, 8, 9\}, \{4, 5, 7, 8, 9\}, \{4, 6, 7, 8, 9\}, \{5, 6, 7, 8, 9\} \}$\\
			\hline
			$d_{322}$ & 6&$\{ \{1, 7\}, \{1, 8\}, \{2, 8\}, \{1, 9\}, \{2, 9\}, \{3, 9\} \}$ \\
			\hline
			$l_4$ &0& $\emptyset$ \\
			\hline
			$l_4\_basis$&0& $\emptyset$ \\
			\hline
			$l_5$ & $3$ & $\{\{2, 3, 4, 5, 6\}, \{3, 4, 5, 6, 7\}, \{4, 5, 6, 7, 8\}\}$\\
			\hline
			$l_5\_basis$& $2$ & $\{\{1, \{2, 3, 4, 5, 6\}\}, \{9, \{4, 5, 6, 7, 8\}\}\}$\\
			 \hline
			 $s_3$ &0& $\emptyset$ \\
			 \hline
			 $s_4$&0& $\emptyset$\\		 
			 \hline
			 $s_5$& $3$ & $\{2, 3, 4, 5, 6\}, \{3, 4, 5, 6, 7\}, \{4, 5, 6, 7, 8\}$\\	 
			 \hline
			 \multirow{5}{*}{$e_3$}& \multirow{5}{*}{$50$}&   $\{\{1, 2, 3\}, \{1, 2, 4\}, \{1, 2, 5\}, \{1, 2, 6\}, \{1, 3, 4\}, \{1, 3, 5\}, \{1, 3, 6\}, \{1, 4, 5\}, \{1, 4, 6\}, \{1, 5, 6\},$\\
			 &&$\{2, 3, 4\}, \{2, 3, 5\}, \{2, 3, 6\}, \{2, 3, 7\}, \{2, 4, 5\}, \{2, 4, 6\},
			 \{2, 4, 7\}, \{2, 5, 6\}, \{2, 5, 7\}, \{2, 6, 7\},$\\
			 &&$ \{3, 4, 5\}, \{3, 4, 6\}, \{3, 4, 7\}, \{3, 4, 8\}, \{3, 5, 6\}, \{3, 5, 7\}, \{3, 5, 8\}, \{3, 6, 7\}, \{3, 6, 8\}, \{3, 7, 8\},$\\
			 && $\{4, 5, 6\}, \{4, 5, 7\}, \{4, 5, 8\}, \{4, 5, 9\}, \{4, 6, 7\}, \{4, 6, 8\}, \{4, 6, 9\}, \{4, 7, 8\}, \{4, 7, 9\}, \{4, 8, 9\},$\\
			 && $  
			 \{5, 6, 7\}, \{5, 6, 8\}, \{5, 6, 9\}, \{5, 7, 8\}, \{5, 7, 9\}, \{5, 8, 9\}, \{6, 7, 8\}, \{6, 7, 9\}, \{6, 8, 9\}, \{7, 8, 9\}\}$\\
			 \hline
			 \multirow{7}{*}{$e_4$}& \multirow{7}{*}{$45$}&  
			$\{\{1, 2, 3, 4\}, \{1, 2, 3, 5\}, \{1, 2, 3, 6\}, \{1, 2, 4, 5\}, \{1, 2, 4, 6\}, \{1, 2, 5, 6\}, 
			\{1, 3, 4, 5\},\{1, 3, 4, 6\}, $\\
			&&$ \{1, 3, 5, 6\}, \{1, 4, 5, 6\}, \{2, 3, 4, 5\}, \{2, 3, 4, 6\},
			\{2, 3, 4, 7\}, \{2, 3, 5, 6\},\{2, 3, 5, 7\}, \{2, 3, 6, 7\},  $\\
			&& $\{2, 4, 5, 6\}, \{2, 4, 5, 7\}, 
			\{2, 4, 6, 7\}, \{2, 5, 6, 7\}, \{3, 4, 5, 6\},\{3, 4, 5, 7\}, \{3, 4, 5, 8\}, \{3, 4, 6, 7\}, $\\
			&&$
			\{3, 4, 6, 8\}, \{3, 4, 7, 8\}, \{3, 5, 6, 7\}, \{3, 5, 6, 8\},\{3, 5, 7, 8\}, \{3, 6, 7, 8\}, 
			\{4, 5, 6, 7\}, \{4, 5, 6, 8\}, $\\
			&&$ \{4, 5, 6, 9\}, \{4, 5, 7, 8\}, \{4, 5, 7, 9\},\{4, 5, 8, 9\},
			\{4, 6, 7, 8\}, \{4, 6, 7, 9\}, \{4, 6, 8, 9\}, \{4, 7, 8, 9\}, $\\
			&&$ \{5, 6, 7, 8\}, \{5, 6, 7, 9\},\{5, 6, 8, 9\}, \{5, 7, 8, 9\}, \{6, 7, 8, 9\}\}$\\
			\hline
			\multirow{3}{*}{$e_5$}& \multirow{3}{*}{$21$}&  $\{\{1, 2, 3, 4, 5\}, \{1, 2, 3, 4, 6\}, \{1, 2, 3, 5, 6\}, \{1, 2, 4, 5, 6\}, \{1, 3, 4, 5, 6\},
			\{2, 3, 4, 5, 6\}, \{2, 3, 4, 5, 7\},$\\
			&&$ \{2, 3, 4, 6, 7\}, \{2, 3, 5, 6, 7\}, \{2, 4, 5, 6, 7\},
			\{3, 4, 5, 6, 7\}, \{3, 4, 5, 6, 8\}, \{3, 4, 5, 7, 8\}, \{3, 4, 6, 7, 8\}, $\\
			&&$\{3, 5, 6, 7, 8\}, \{4, 5, 6, 7, 8\}, \{4, 5, 6, 7, 9\}, \{4, 5, 6, 8, 9\}, \{4, 5, 7, 8, 9\}, \{4, 6, 7, 8, 9\},\{5, 6, 7, 8, 9\}\}$\\
			\hline
			$i_2$ &0& $\emptyset$\\
			\hline
			$se_6$ &$3$& $\{\{1, 2, 3, 4, 5, 6\}, \{2, 3, 4, 5, 6, 7\}, \{3, 4, 5, 6, 7, 8\}, \{4, 5, 6, 7, 8, 9\}\}$\\
			\hline
			$se_7$ &0& $\emptyset$\\
			\hline
			\multirow{12}{*}{$g_{322}$}& \multirow{12}{*}{12}& \multicolumn{1}{c|}{$(1 2 3 4 5 6 7 8 9)$}\\
			&& \multicolumn{1}{c|}{$(1 2 3 4 6 5 7 8 9)$}\\
			&& \multicolumn{1}{c|}{$(1 2 3 5 4 6 7 8 9)$}\\
			&& \multicolumn{1}{c|}{$(1 2 3 5 6 4 7 8 9)$}\\
			&& \multicolumn{1}{c|}{$(1 2 3 6 4 5 7 8 9)$}\\
			&& \multicolumn{1}{c|}{$(1 2 3 6 5 4 7 8 9)$}\\
			&& \multicolumn{1}{c|}{$(9 8 7 4 5 6 3 2 1)$}\\
			&& \multicolumn{1}{c|}{$(9 8 7 4 6 5 3 2 1)$}\\
			&& \multicolumn{1}{c|}{$(9 8 7 5 4 6 3 2 1)$}\\
			&& \multicolumn{1}{c|}{$(9 8 7 5 6 4 3 2 1)$}\\
			&& \multicolumn{1}{c|}{$(9 8 7 6 4 5 3 2 1)$}\\
			&& \multicolumn{1}{c|}{$(9 8 7 6 5 4 3 2 1)$}\\
				\Xcline{1-3}{1.2pt}
			
		\end{tabular}
	}
	
  	\caption{Combinatorics of $P_{322}$.}
  	\label{table:p322}
\end{table}

\begin{figure}[H]
  	\scalebox{0.4}[0.4]{\includegraphics {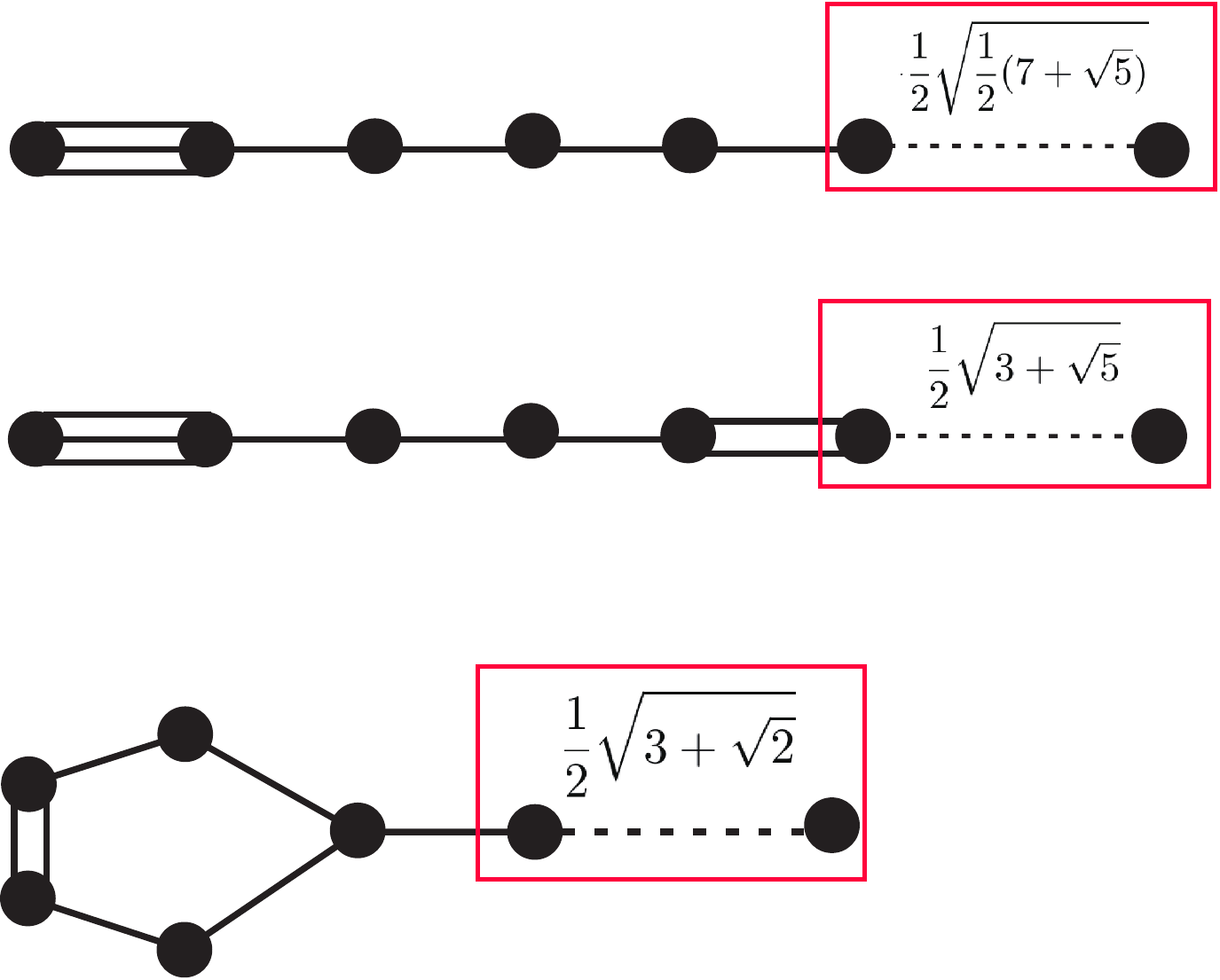}}
  	\caption{Compact simplicial prisms in $\mathbb{H}^5$.}
  	\label{figure:prism5}
  \end{figure}

It is worthy to mention that the set $l_5\_basis$, if not empty, can help to reduce the computation since the list of simplicial $5$-prisms is available. For example, suppose $\{2,4,5,6,7\}$ (referring to facets $F_2,F_4,F_5,F_6,F_7$) bound a facet $8$ (means $F_8$) of 4-simplex. Then we can assume that $F_8$ is orthogonal to $F_2$, $F_4$, $F_5$, $F_6$, and $F_7$ (i.e., $m_{28}=m_{48}=m_{58}=m_{68}=m_{78}=2$). The vectors obtained this way can be treated as \emph{bases}, named basis vectors, and all of the other potential vectors that may lead to a Gram vector can be realized by gluing the simplicial $5$-prisms, as shown in Figure \ref{figure:prism5}, at their orthogonal ends.


Moreover, among all of the $322$ polytopes, we only need to study those with number of hyperparall pairs larger than or equal to two due to the following theorems:

\begin{theorem}[\cite{FT:08},part of Theorem A]
	Every compact hyperbolic Coxeter $d$-polytope, $d>4$, has a pair of disjoint facets.
\end{theorem}

\begin{theorem}[\cite{FT:09},Main Theorem A]
	A compact hyperbolic Coxeter $d$-polytope with exactly one pair of non-intersecting facets has at most $d + 3$ facets.
\end{theorem}

There are $109$ polytopes with two or more hyperparallel facets. We group them by the number $d_k$ of disjoint pairs as illustrated in Table \ref{table: group}.


 \begin{table}[h]
 	{\footnotesize
 	\begin{tabular}{c|c}
 		\Xcline{1-2}{1.2pt}
 		\textbf{$d_k$} & \textbf{labels of polytopes}\\
 		\hline
 		$6$ & 312 319 322 \\
 		\hline
 		$5$ &254 302 311 315 318 320 321 \\
 		\hline
 		$4$ & 249 255 280 291 292 293 295 300 301 303 310 313 314 316 317\\
 		\hline
 			\multirow{2}{*}{3} &2 10 175 240 241 242 246 247 248 252 256 264 268 271 272 273 279 282 284 286 287 289 290 294 296\\
 		&
 		297 298 299 304 305 306 307 308 309 \\
 		\hline
 			\multirow{2}{*}{2} & 4 9 18 61 68 110 173 174 176 179 180 182 183 203 209 210 212 214 216 218 225 234 239 243 244 245 \\
 			&
 		250 251 253 257 258 259 260 261 262 263 265 266 267 269 270 274 275 276 277 278 281 283 285 288\\
 		\Xcline{1-2}{1.2pt}
 		
 	\end{tabular}
 	}
 	
 	\hspace*{0.5cm}
 	\caption{Five groups with respect to different numbers of disjoint pairs.}
 	\label{table: group}
 \end{table}

\section{Block-pasting algorithms for enumerating all of the candidate matrices over a certain combinatorial type.} \label{chapter:algorithm}

We now use the \emph{block-pasting algorithm} to determine all of the potential matrices for  the $109$ compact  combinatorial types reported in Section \ref{section:4d8f}. Recall that the entries have only finite options, i.e., $k_{ij}\in\{1,2,3,\cdots, 7\}\cup \{\infty\}$. Compared to the backtracking search algorithm raised in \cite{JT:2018}, ``block-pasting" algorithm is more efficient and universal. Generally speaking, the backtracking search algorithm uses the method of ``a series circuit", where the potential matrices are produced one by one. Whereas, the block-pasting algorithm adopts the idea of ``a parallel circuit", where different parts of a potential matrix are generated simultaneously and then pasted together. 

For each  vertex $v_i$ of a $5$-dimensional hyperbolic Coxeter polytope $P_k$, we define the \emph{i-chunk}, denoted as $k_i$, to be the ordered set of the five facets intersecting at the vertex $v_i$ with increasing subscripts.  For example, for the polytope $P_{322}$ discussed above, there are $18$ chunks as it has $18$ vertices. We may also use $k_i$ to denote the ordered set of subscripts, i.e., $k_i$ is referred to as a set of integers of length five.

Since the compact hyperbolic $5$-dimensional polytopes are simple, each chunk possesses $\tbinom{5}{2}=10$ dihedral angles, namely the angles between every two adjacent facets. For every chunk $k_i$, we define an \emph{i-label} set $e_i$ to be the ordered set $\{10a+b|\{a,b\}\in E_i\}$, where $E_i$, named the \emph{i-index} set, is the ordered set of pairs of facet labels. These are formed by  choosing every two members from the chunk $k_i$ where the labels increase lexicographically. For example, suppose the five facets intersecting at the first vertex are $F_1$, $F_2$, $F_3$, $F_5$, and $F_6$. Then, we have 
$$k_1=\{F_1,F_2,F_3,F_5,F_6\} (\text{or} \{1,2,3,5,6\}),$$ $$E_1=
\{\{1,2\},\{1,3\},\{1,5\},\{1,6\},\{2,3\},\{2,5\},\{2,6\},\{3,5\},\{3,6\},\{5,6\}\},$$ $$e_1=
\{12,13,15,16,23,25,26,35,36,56\}.$$

Next, we list all of the Coxeter vectors of the elliptic Coxeter diagrams of rank $5$. Note that we have made the convention of considering only the diagrams with integer entries less than or equal to seven. The qualified Coxeter diagrams are shown in Figure \ref{figure:elliptic5}:

\begin{figure}[h]
	\scalebox{0.3}[0.3]{\includegraphics {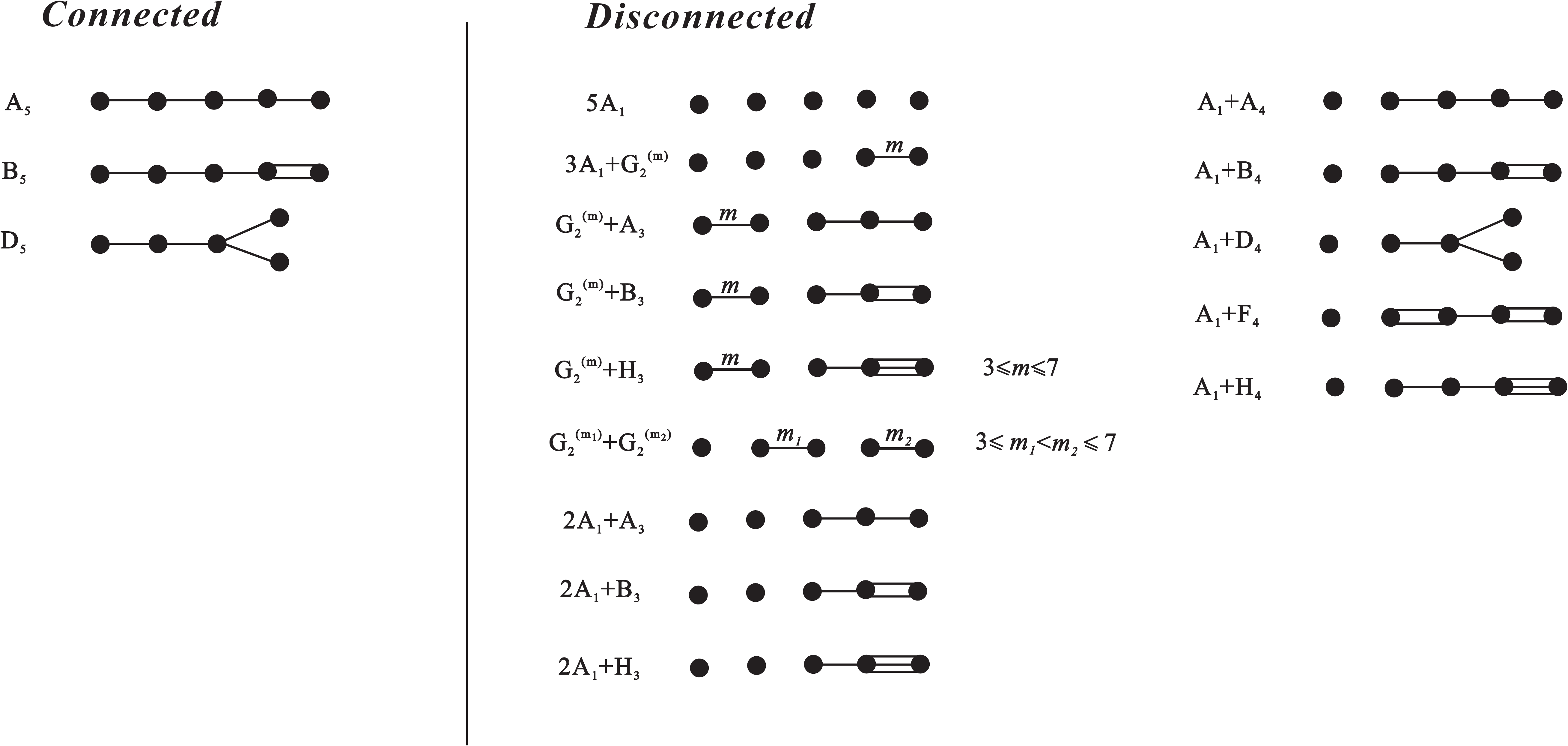}}
	\caption{Spherical Coxeter diagram of rank $5$ with labels less than or equal to seven.}\label{figure:elliptic5}
\end{figure}

\begin{figure}[h]
	\scalebox{0.4}[0.4]{\includegraphics {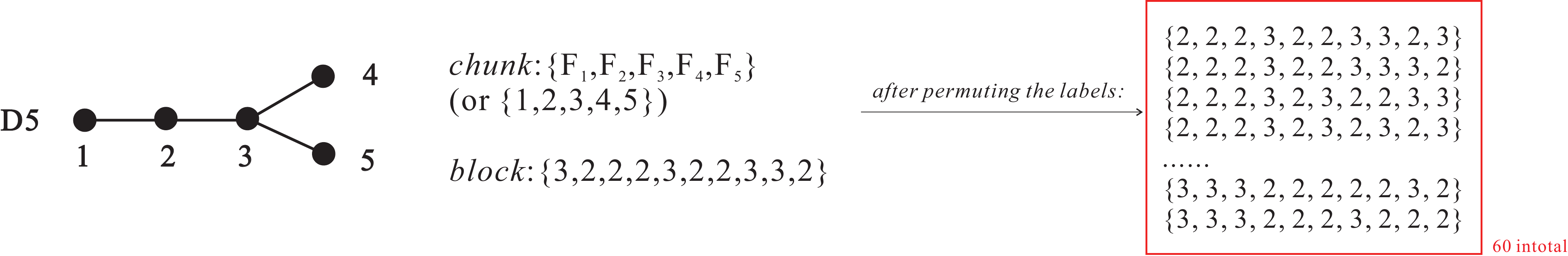}}
	\caption{Building the pre-block $\mathcal{S}$.}\label{figure:d5}
\end{figure}

We apply the permutation group on five letters $S_5$ to the labels of the nodes of the Coxeter diagrams in Figure \ref{figure:elliptic5}. This produces all of the possible vectors when varying the order of the five facets. For example, there are $60$ vectors for the single diagram $D_5$ as shown in Figure \ref{figure:d5}.  There are $1946$ distinct such vectors of rank 5 elliptic Coxeter diagrams. The set of all of these vectors is called the \emph{pre-block}; it is denoted by $\mathcal{S}$. The set $\mathcal{S}$ can be regarded as a $1946\times 10$ matrix in the obvious way. In the following, we do not distinguish these two viewpoints and may refer to $\mathcal{S}$ as either a set or a matrix. 

 We then generate every datafram $B_i$, named the $i$-block, of size $1946\times 10$, corresponding to every chunk $k_i$, for a given polytope $P_k$, where $1\leq i \leq |V_k|$ and $|V_k|$ is the number of vertices of $P_k$. Firstly we evaluate $B_i$ by $\mathcal{S}$ and take the ordered set $e_i$ defined above as the column names of $B_i$. For example, for $e_1=\{12,13,15,16,23,25,26,35,36,56\}$, the columns of $B_i$ are referred to as $(12)$-, $(13)$-, $(15)$-, $(16)$-, $(23)$-, $(25)$-, $(26)$-, $(35)$-, $(36)$-, $(56)$-columns. 
 
 Denote $L$ to be a vector of length $36$ as follows:
$$L=\{12,13,...,19{\color{red},}~23,24,...,29{\color{red},}~34,35,...,39{\color{red},}~45,46,...,49{\color{red},}~56,57,...,59{\color{red},}~67,68,69{\color{red},}~ 78,79{\color{red},}~89\}.$$
Then all of the numbers in the label set of $d_k$ (the set of disjoint pair of facets of $P_k$) are excluded from $L$ to obtain a new vector For brevity, the new vector is also denoted by $L$. For example, the numbers excluded for the polytope $P_{322}$ are $17$, $18$, $28$, $19$, $29$, and $39$ as illustrated in Table \ref{table:p322}. The length of $L$ is denoted by $l$.  

Next, we extend every $1946\times 10$ dataframe to a $1946\times l$ one, with column names $L$, by simply putting each $(ij)$-column to the position of corresponding labelled column, and filling in the value zero in the other positions. We continue to use the same notation $B_i$ for the extended dataframe. In the rest of the paper, we always mean the extended dataframe when using the notation $B_i$. 

After preparing all of the the blocks $B_i$ for a given polytope $P_k$, we proceed to paste them up. More precisely, when pasting $B_1$ and $B_2$, a row from $B_1$ is matched up with a row of $B_2$ where every two entries specified by the same index $i$, where $i\in e_1\cap e_2$, have the same values. The index set $e_1\cap e_2$ is called a \emph{linking key} for the pasting. The resulting new row is actually the sum of these two rows in the non-key positions; the values are retained in the key positions. The dataframe of the new data is denoted by $B_1\cup^*B_2$.

We use the following example to explain this process. Suppose

$B_1=\{x_1,x_2\}=\{(1,2,4,4,2,6,0,0,\cdots ,0), (1,2,4,5,2,6,0,0,\cdots,0)\}$,

$B_2=\{y_1,y_2,y_3\}=\{(1,2,4,4,0,0,1,7,0,0,\cdots,0),(1,2,4,4,0,0,6,5,0,0,\cdots,0)$,

\hspace{3.5cm}$(1,2,3,4,0,0,1,7,0,0,\cdots,0)\} .$

In this example, $x_1$ has the same values with $y_1$ and $y_2$ on the $(12)$-, $(13)$-, $(14)$- and $(15)$- positions. In other words, the linking key here is $\{12,13,14,15\}$. Thus, $y_1$ and $y_2$ can paste to $x_1$, forming the Coxeter vectors $$(1,2,4,4,2,6,1,7,0,...,0)~\text{and}~(1,2,4,4,2,6,6,5,0,...,0), \text{respectively}.$$
\noindent In contrast, $x_2$ cannot be pasted to any element of $B_2$ as there are no vectors with entry $5$ on the $(15)$-position. Therefore, $$B_1\cup^* B_2=\{(1,2,4,4,2,6,1,7,0,0,\cdots,0),(1,2,4,4,2,6,6,5,0,0,\cdots,0)\}.$$

We then move on to paste the sets $B_1 \cup^*B_2$ and $B_3$. We follow the same procedure with an updated index set. Namely the linking key, is now $e_1\cup e_2\cap e_3$. We conduct this procedure until we finish pasting the final set $B_{\vert V_k\vert}$. The set of linking keys used in this procedure is $$\{e_1\cap e_2, e_1\cup e_2\cap e_3,\cdots,e_1\cup e_2\cup\cdots\cup e_{i-1}\cap e_i,\cdots ,e_1\cup e_2\cup \cdots \cup e_{\vert V\vert -1}\cap e_{\vert V\vert}\}.$$

After pasting the final block $B_{\vert V_k\vert}$, we obtain all of the potential vectors of the given polytope. This approach has been Python-programmed on a PARATERA server cluster.

 When carrying out this approach, all of the potential vectors of $P_{322}$ are enumerated successfully. However it encounters a memory error in some cases. For example, for the polytope $P_{10}$, the computation incurs an error when pasting $B_{16}$ as illustrated in Table \ref{table:p10}.
 
 \begin{table}[h]
 	{\footnotesize
 		\begin{tabular}{c|c|c|c}
 			\Xcline{1-4}{1.2pt}
 			\multicolumn{2}{c|}{\quad \quad \textbf{$i$-chunk ($1\leq i\leq 24$)\quad \quad }}
 			
 			& \quad \quad 	\textbf{data size after pasting $b_i$} \quad \quad \quad \quad  & \textbf{time consumed (s)} \quad \quad \\ \hline
 		$h_1$& $12345$ & $1946$ &$0$\\	\hline
 		$h_2$& $12346$ & $32780$ &$0.05$\\	\hline
 		$h_3$& $12356$ & $46851$ &$0.05$\\	\hline
 		$h_4$& $12457$ & $909145$ &$0.7$\\	\hline
 		$h_5$& $12467$ & $1295489$ &$1.06$\\	\hline
 		$h_6$& $12567$ & $774519$ &$0.69$\\	\hline
 		$h_7$& $13457$ & $1619039$ &$1.32$\\	\hline
 		$h_8$& $13467$ & $1137673$ &$1.06$\\	\hline
 		$h_9$& $13567$ & $1052454$ &$0.79$\\	\hline
 		$h_{10}$& $23458$ & $20832842$ &$16.68$\\	\hline
 		$h_{11}$& $23468$ & $34761662$ &$24.93$\\	\hline
 		$h_{12}$& $23568$ & $20789435$ &$17.63$\\	\hline
 		$h_{13}$& $24578$ & $40260708$ &$28.12$\\	\hline
 		$h_{14}$& $24678$ & $29788033$ &$21.62$\\	\hline
 		$h_{15}$& $25678$ & $27886613$ &$16.95$\\	\hline
 		$\cdots$ & $\dots$& Memory error &$\dots$\\
 		\Xcline{1-4}{1.2pt}

 		\end{tabular}
 }
 
 \hspace*{0.5cm}
 	\caption{Pasting of $P_{10}$.}\label{table:p10}
 \end{table}
 

 Therefore, a refined algorithm is needed to continue this research. The philosophy of the refinement is to introduce more necessary conditions other than the vertex spherical and Euclidean square restrictions, to reduce the number of vectors during the block-pasting. Firstly, we collect data sets $\mathcal{L}_4$, $\mathcal{L}_5$, $\mathcal{L}_5\_basis$, $\mathcal{S}_3$, $\mathcal{S}_4$, $\mathcal{S}_5$, $\mathcal{S}_6$, $\mathcal{S}_7$, $\mathcal{E}_3$, $\mathcal{E}_4$,
 $\mathcal{E}_5$, $\mathcal{E}_6$, $\mathcal{E}_7$, $\mathcal{I}_2$, as claimed in Table \ref{table:library}, by the following two steps:
\begin{enumerate}
	\item Prepare Coxeter diagrams of rank $r$, as assigned in Table \ref{table:library}, and write down the Coxeter vectors under an arbitrary system of node labeling.
	\item Apply the permutation group $S_r$ to the labels of the nodes and produce the desired data set consisting of all of the distinct Coxeter vectors under all of the different labelling systems.
\end{enumerate}

Note that the set $\mathcal{S}_5$ is exactly the pre-block $\mathcal{S}$ we construct before. Readers can refer to the process of producing $\mathcal{S}$ for the details of building these data sets.
  
\begin{table}[H]
	{\footnotesize
			\begin{tabular}{c|c|c|c}
			\Xcline{1-4}{1.2pt}
			\multirow{2}{*}{\textbf{Types of Coxeter diagrams}}
			 & 	\textbf{\# Coxeter}  & \textbf{\# distinct Coxeter Vectors } 
			& \multirow{2}{*}{\textbf{data sets}}\\
			&\textbf{diagrams}& \textbf{after permutation on nodes} & \\   
			\hline	
			Coxeter diagrams of
			compact hyperbolic 4-simplex& 5 & 420 &$\mathcal{L}_5$\\
			\hline
			Coxeter diagrams of
			compact hyperbolic 3-simplex& 9& 108 &$\mathcal{L}_4$\\
			\hline
			rank $3$ elliptic Coxeter diagrams & 9& 31&$\mathcal{S}_3$ \\
			\hline
			rank $4$ elliptic Coxeter diagrams & 29& 242 & $\mathcal{S}_4$ \\
			\hline
			rank $5$ elliptic Coxeter diagrams &47 & 1946 &$\mathcal{S}_5$\\
			\hline
			rank $6$ elliptic Coxeter diagrams & 117& 20206 &$\mathcal{S}_6$\\
			\hline
			rank $7$ elliptic Coxeter diagrams & 196 & 227676 & $\mathcal{S}_7$  \\
			\hline
			rank $3$ connected parabolic Coxeter diagrams &3 & 10 &$\mathcal{E}_3$ \\
			\hline
		    rank $4$ connected parabolic Coxeter diagrams & 3 & 27 &$\mathcal{E}_4$\\	
			\hline
			rank $5$ connected parabolic Coxeter diagrams&5&257&
			$\mathcal{E}_5$\\
			\hline
			rank $6$ connected parabolic Coxeter diagrams&4&870& $\mathcal{E}_6$\\
			\hline
			rank $7$ connected parabolic Coxeter diagrams &5&6870& $\mathcal{E}_7$\\
			\hline
			Coxeter diagrams of Euclidean $2$-cube &4&3&$\mathcal{I}_2$\\
			\Xcline{1-4}{1.2pt}
			\end{tabular}
	}

		\caption{Data sets used to reduce the computational load.}
		\label{table:library}
	\end{table}


Next, we modify the block-pasting algorithm by using additional metric restrictions. More precisely, remarks \ref{remark:1}--\ref{remark:4}, which are practically reformulated from Theorem \ref{Vinberg:thm3.1}, must be satisfied. 

\begin{remark}(``$l4$/$l5$-condition") \label{remark:1}
	The Coxeter vector of the six/ten dihedral angles formed by the four/five facets with the labels indicated by the data in $l_4$/$l_5$ is {\color{red} IN} $\mathcal{L}_4$/$\mathcal{L}_5$
\end{remark}

\begin{remark}(``$s3$/$s4$/$s5$/$s6$/$s7$-condition") \label{remark:2}
	The Coxeter vector of the three/six/ten/fifteen/twenty-one dihedral angles formed by the three/four/five/six/seven facets with the labels indicated by the data in $s_3$/$s_4$/$s_5$/$se_6$/$se_7$ is {\color{red} NOT IN} $\mathcal{S}_3$/$\mathcal{S}_4$/$\mathcal{S}_5$/$\mathcal{S}_6$/$\mathcal{S}_7$.
\end{remark}

\begin{remark}(``$e3$/$e4$/$e5$/$e6$/$e7$-condition") \label{remark:3}
	The Coxeter vector of the three/six/ten/fifteen/twenty-one dihedral angles formed by the three/four/five/six/seven facets with the labels indicated by the data in $e_3$/$e_4$/$e_5$/$se_6$/$se_7$ is {\color{red} NOT IN} $\mathcal{E}_3$/$\mathcal{E}_4$/$\mathcal{E}_5$/$\mathcal{E}_6$/$\mathcal{E}_7$.
\end{remark}

\begin{remark}(``$i2$-condition") \label{remark:4}
	The Coxeter vector formed by the four facets with the labels indicated by the data in $i_2$ is {\color{red} NOT IN} $\mathcal{I}_2$.
\end{remark}

The {\color{red}``IN"} and {\color{red}``NOT IN"} tests are called the ``saving" and the ``killing" conditions, respectively. The ``saving" condition is much more efficient than the ``killing" condition because the ``what kinds of vectors are qualified" is much more restrictive than the ``what kinds of vectors are not qualified". The more delicate ``saving conditions" are introduced in next section. Moreover, we remark that the $l3$-condition and the sets $l_3$ and $\mathcal{L}_3$, which can be defined analogously as the $l4$ / $l5$ setting, are not introduced. This is because the effect of using both $s3$- and $e3$- conditions is equivalent to adopting the $l3$-condition.   

We now program these conditions and insert them into the appropriate layers of the pasting to reduce the computational load. Here the ``appropriate layer" means the layer where the dihedral angles are non-zero for the first time. For example, for $\{3,6,7\}\in e_3$, we find the $j$-th block pasting where the data in columns ($36$-,$37$-,$67$-) of the dataframe $B_1\cup^*B_2\cdots \cup^*B_j$ become non-zero. Therefore, the $e3$-condition for $\{3,6,7\}$ is inserted immediately after the $j$-th block pasting. The symmetries of the polytopes are factored out when the pastes are finished. The matrices (or vectors) after all these conditions (metric restrictions and symmetry equivalence) are called  \emph{``SEILper"-potential matrices (or vectors)} of certain combinatorial types. All of the numbers of the results are reported in Tables \ref{table:g6}--\ref{table:g2}. In those tables, the notation $grp_i(j,k)$ means there are $j$ polytopes with $i$ pairs of disjoint facets, and $k$ of them admit SEILper potential vectors, i.e., the number of SEILper potential vectors that are non-zero. The letters ``$s$, $m$, $h$, and $d$" on the lines named ``time" represent ``seconds, minutes, hours, and days", respectively.  The numbers in red indicate that the corresponding polytopes have a non-empty set $l_5\_basis$; therefore, the results obtained are basis SEILper potential vectors. This calculation is called the \emph{basis approach}. We confirm these cases without using the $l5$-condition, called the \emph{direct approach}, in the validation part presented in Section \ref{section:vadilation}. 
  
\begin{table}[h]
	{\footnotesize
		\begin{tabular}{c|ccc}

			\Xcline{1-4}{1.2pt}	
			grp6 (3/3)&	{\color{red}312 }& {\color{red}319} & {\color{red} 322}\\
			\hline
			\#SEILper&	1&	5&	3\\
			\hline
			time&	37.755 s&	26.48 s&	23.495 s\\
			\Xcline{1-4}{1.2pt}

		\end{tabular}
	}
	
	\hspace*{0.5cm}
	\caption{Result of group 6. }
	\label{table:g6}
\end{table}

\begin{table}[h]
	{\footnotesize
		\begin{tabular}{c|ccccccc}
			\Xcline{1-8}{1.2pt}
			
			grp5 (7/7) &	254&	{\color{red}302}&	{\color{red}311}&	{\color{red}315}	& {\color{red}318}&	{\color{red}320}	& {\color{red}321}\\
			\hline
			\#SEILper&	9467&	5&	288&	49&	2&	6&	3\\
			\hline
			time&	2.925 h	&45.49 s&	3.405 m&	42.06 s&	33.81 s&	53.115 s&	8.685 m\\
			
			\Xcline{1-8}{1.2pt}

		\end{tabular}
	}
	
	\hspace*{0.5cm}
 	\caption{Result of group 5. }
	\label{table:g5}
\end{table}

  \begin{table}[h]
  	{\footnotesize
  		\begin{tabular}{c|ccccccccc}
  			\Xcline{1-10}{1.2pt}
  			
  		gr4 (15/9)	&249&	255&	280&	291&	292	&293&	{\color{red}295}&	{\color{red}300}&{\color{red}	301}\\
  		\hline
  		\#SEILper &3910&	282&	148	&37&	10&	0&	0&	383&	5\\
  		\hline
  		time&	5.49 m	&5.6 m&	1.46 h&	2.03 m&	8.91 m&	4.2 m&	5.11 m&	9.755 h&	2.675 m\\
  		\Xcline{1-10}{1.2pt}
  		\multicolumn{1}{c}{}&{\color{red}303}&	{\color{red}310}	& {\color{red}313}&	{\color{red}314	}& {\color{red}316}&	{\color{red}317}&&&\\
  		\cline{2-10}
  			\multicolumn{1}{c}{}&0&	0&	68&	0&	8&	0&&&\\
  		\cline{2-10}
  			\multicolumn{1}{c}{}&5.23 m	&5.64 m&	8.215 m	&1.805 m&	0.855 m&	8.5 m&&&\\
  		
  		\Xcline{1-10}{1.2pt}

  		\end{tabular}
  	}
  	
  	\hspace*{0.5cm}
  	\caption{Result of group 4. }
  	\label{table:g4}
  \end{table}  

    \begin{table}[h]
    	{\footnotesize
    		\begin{tabular}{cccccccccc}
    			\Xcline{1-10}{1.2pt}
    			
    			\multicolumn{1}{c|}{grp3 (34/16)}&	2&	10&	175&	240	&241&	242&	246&	247&	248\\
    			\hline
    			\multicolumn{1}{c|}{\#SEILper}	&	160395&	139143	&0&	4&	54&	6744&	0&	1546&	2574\\
    			\hline

    			\multicolumn{1}{c|}{time}&	7.545 d&	15.585 h&	10.325 s&	17.815 m&	6.07 m&	10.015 h&	4.8 m&	8.63 m&	15.845 m\\
    			
    			\hline
   
    			\hline
    			&252	&256&	264	&268&	271&	272	&273&	279&	282	\\
    			\cline{2-10}
    			&0&	15&	514&	155&	0&	0&	0&	0&	40\\
    			\cline{2-10}
    		&	3.6 m&	3.615 m&	10.29 h&	0.63 h&	4.26 m&	7.455 s&	3.155 m&	7.735 m&	20.96 m\\
    			 
    		\Xcline{2-10}{1.2pt}
    			     			
    		&284&	286	&287&	289&	290&	{\color{red}294}&	{\color{red}296}&	{\color{red}297}&	{\color{red}298}\\
    		\cline{2-10}
    		&60	&0	&0&	0&	0&	83&	37&	0&	11\\
    		\cline{2-10}
    		&4.145 m&	3.39 m&	2.665 m&	6.655 m&	7.46 m&	4.975 m&	0.96 d&	17.325 m&	7.165 h\\
    		\Xcline{2-10}{1.2pt}
    		&{\color{red}299}&	{\color{red}304}	& {\color{red}305}&	{\color{red}306}&	{\color{red}307}&	{\color{red}308}&	{\color{red}309}&&\\
    		\cline{2-10}
    	&	0&	0&	0&	0&	0&	96&	0&&\\
    	\cline{2-10}
    		&1.975 h&	3.195 m&	2.685 h&	2.455 h&	2.125 h&	5.89 h	&3.185 h&&\\
    		\Xcline{1-10}{1.2pt}

    		\end{tabular}
    	}
    	
    	\caption{Result of group 3. }
    	\label{table:g3}
    \end{table}

 \begin{table}[h]
 	{\footnotesize
 		\begin{tabular}{cccccccccc}
 			\Xcline{1-10}{1.2pt}
			\multicolumn{1}{c|}{grp2(50/15)}  &	4&	9&	18&	61&	68&	110	&173&	174&	176\\
			\hline
			\multicolumn{1}{c|}{\#SEILper}&	11420&	491	&1058&	79&	805&	0&	0	&140&	0\\
			\hline
			\multicolumn{1}{c|}{time}&	7.885 h&	2.025 d&	13.11 d&	8.27 d&	3.33 d&	2.68 h&	2.56 h&	8.025 h	&20.305 s\\
			\Xcline{1-10}{1.2pt}
			&179&	180	&182&	183&	203	&209&	210&	212&	214\\
			\cline{2-10}
			&0&	1371&	2336&	0&	908&	0&	0	&0	&2504\\
			\cline{2-10}
			&3.31 m&	8.255 m	&24.205 d&	2.99 m&	19.885 d&	5.565 m&	5.965 m&	23.29 s&	11.665 m\\
			\Xcline{2-10}{1.2pt}
			&216&	218	&225&	234	&239&	243&	244&	245&	250\\
			\cline{2-10}
			&260&	40&	0&	0&	0&	0&	0&	0&	0\\
			\cline{2-10}
			&12.545 m&	12.2 m&	3.095 m&	1.635 m	&2.61 m&	4.785 h&	6.74 m&	1.04 m&	7.595 h\\
			\Xcline{2-10}{1.2pt}
			&251&	253	&257&	258&	259	&260&	261	&262&	263\\
			\cline{2-10}
			&0&	0&	25&	70&	0&	0&	0&	0&	0\\
			\cline{2-10}
			&1.15 h&	3.755 m	&13.4 m&	16.075 h&	3.265 h	&3.835 m&	0.835 h	&5.22 h&	0.845 h\\
			\Xcline{2-10}{1.2pt}
			&265&	266&	267&	269&	270&	274&	275&	276&	277\\
			\cline{2-10}
			&0&	0&	0&	0&	0&	0&	0&	96&	0\\
			\cline{2-10}
			&2.655 m&	13.97 s&	3.965 m&	2.855 h&	2.96 h&	4.495 h&	3.165 h	&10.155 h&	5.02 h\\
			\Xcline{2-10}{1.2pt}
			&278&	281&	283	&285&	288	&&&&\\
			\cline{2-10}			
			&0	&0&	0&	0&	0	&&&&\\
			\cline{2-10}			
			&1.17 h	&1.28 h&	4.485 m	&1.91 h&	3.17 h&&&&\\				
            \Xcline{1-10}{1.2pt}
			\cline{1-9}
			\end{tabular}
			}
			
			\caption{Result of group 2. }
			\label{table:g2}
			\end{table}

So far, we have restricted the number of candidates and perform the computations in a reasonable amount of time. However, when passing to the second phase of calculating the signature, there are cases, such as the  $9467$ SEILper potential matrices of the polytoep $P_{254}$, that require an unsatisfying amount of time. In the next section,  the restriction of ``intersection types" is introduced to facilitate the computation of signature in $30s$ per polytope.

\section{Signature Constraints of hyperbolic Coxeter \texorpdfstring{$n$}-polytopes}\label{section:signature}

After preparing all of the SEILper matrices, we proceed to calculate the signatures of the potential Coxeter  matrices to determine if they lead to the Gram matrix $G$ of an actual hyperbolic Coxeter polytope. 

Firstly, we modify every SEILper matrix $M$ as follows: 

\begin{enumerate}
	\item Replace $\infty$s by unknowns of $x_i$;
	\item Replace $2$, $3$, $4$, $5$, and $6$ by $0$, $-\frac{1}{2}$, $-\frac{l}{2}$, $-\frac{m}{2}$, $-\frac{n}{2}$, where $$l^2-1=2,~l>0,~m^2-m-1=0,~m>0,~n^2-3=0,~n>0;$$
	\item Replace $7$s by unknowns of $\frac{y_i}{2}$.
\end{enumerate} 

By Theorem \ref{thm:signature}, the resulting Gram matrix must have signature $(5, 1)$. This implies that the determinant of every $6\times 6$ minor of each modified $9\times 9$ SEILper matrix is zero. Therefore, we have the following system of equations and inequality on $x_i$, $l$, $m$, $n$, and $y_i$ to further restrict and lead to the Gram matrices of the desired polytopes:
$$(6.1) ~~~
\begin{cases}
2\det (M_i)=0,\text{for~ any of the~}\tbinom{9}{2}=36 ~6\times 6~\text{minor}~M_i~\text{of}~ M. \\
1.8<y_i<2 ~\text{for~ all}~y_i\\
x_i<-1~\text{for~ all}~x_i\\
l^2-1=2,~l>0,~m^2-m-1=0,~m>0,~n^2-3=0,~n>0
\end{cases}$$

The above conditions are initially stated by Jacquemet and Tschantz in \cite{JT:2018}. Due to practical constraints in \emph{Mathematica}, we denote $2cos(\frac{\pi}{4}),~2cos(\frac{\pi}{5}),~2cos(\frac{\pi}{6})$ by $\frac{l}{2},~\frac{m}{2},~\frac{n}{2}$, rather than $l,~m,~n$ and set $2\det(M_i)=0$ rather than $\det(M_i)=0$. Delicate reasons for doing so can be found in \cite{JT:2018}. Moreover, we first find the \emph{Gr\"{o}bner bases} of the polynomials involved, i.e., $2\det(M_i),~l^2-1,~m^2-m-1,~ n^2-3$, before solving the system. This quickly passes over the cases that have no solution. However, when dealing with some combinatorial types, like $P_{254},~ P_{249},~P_{2},~P_{10}$, etc., the computation cannot be accomplished in reasonable amount of time. In some cases, a single matrix can require more than two hours to compute, which is costly and impedes the validation process. Hence, we introduce the "\emph{intersection restrictions}"  to further reduce the number of SEILper matrices.

Let $P^3$, $P^4$, and $D^4$ denote the combinatorial types of simplicial $3$-prism, simplicial $4$-prism, and the product of two $2$-simplices. We prepare some data for these three combinatorial type as shown in Table \ref{table:p3p4d4}.

\begin{table}[h]
	{\footnotesize
		\begin{tabular}{c|c|c|c}
			\Xcline{1-4}{1.2pt}
			Combinatorial type &$P^3$ & $P^4$ & $D^4$\\
			\hline
			\# facets & 5 &6 &6\\
			\hline
			\multirow{1}{*}{datum}& $s_3$ $s_4$ &  $s_3$ $s_4$ $s_5$  &  $s_3$ $s_4$  $s_5$ \\
			($d_*$:set of disjoint pair of facets;&$e_3$ $e_4$  & $e_3$ $e_4$ $e_5$ & $e_3$ $e_4$ $e_5$\\
			$g_*$: symmetry group )&$i_2$ $d_{P^3}$ $g_{P^3}$ & $i_2$ $l_4$ $l_4\_basis$ $d_{P^4}$ $g_{P^4}$& $i_2$ $l_4$ $l_4\_basis$ $d_{D^4}$ $g_{D^4}$\\
			\Xcline{1-4}{1.2pt}
		\end{tabular}
		}
		
       \hspace*{0.5cm}
       \caption{Combinatorics  of $P^3$, $P^4$, and $D^4$ that are used to obtain SEILper matrices. }
	   \label{table:p3p4d4}
  \end{table}
  
We use $\mathcal{S}_3$, $\mathcal{S}_4$, and $\mathcal{S}_4$ as the pre-blocks of $P^3$, $P^4$ and $D^4$, and apply remarks \ref{remark:1}--\ref{remark:4} to the data presented in Table \ref{table:p3p4d4} to obtain the SEILper potential matrices. The calculation results are arranged on the first line of Table \ref{table:sig344}. 

Next, we generate the sets $\mathcal{P}_3$, $\mathcal{P}_4$, and $\mathcal{D}_4$ as defined below:

\begin{enumerate}
	\item  The set $\mathcal{P}_3$ consists of all of the SEILper potential matrices of the simplicial $3$-prism $P^3$ with the signature $(3,1)$ or $(4,1)$. Equivalently, we are finding the Coxeter groups corresponding to $P^3$ that admit the Fuchsian or quasi-Fuchsian representations into $PO_{4,1}(\mathbb{R})$.
	\item The set $\mathcal{P}_4$ consists of all of the SEILper potential matrices of the simplicial $4$-prism $P^4$ with the signature $(4,1)$ or $(5,1)$.
	\item The set $\mathcal{D}_4$ consists of all of the SEILper potential matrices of the product of two $2$-simplices $D^4$ with the  signature $(4,1)$ or $(5,1)$.
\end{enumerate}

We now explain the signature calculation method for $D^4$, and it is analogous for $P^3$ and $P^4$. From the combinatorial type $D^4$, the signature of the potential matrix has at least four positive eigenvalues. Therefore the signature is one of $(4,0)$, $(4,1)$,$(5,0)$, $(4,2)$, $(5,1)$, and $(6,0)$. We modify the SEILper matrices as described in the beginning of this section and suppose the eigenvalue of the given matrix to be $r_1,r_2,\cdots, r_6$. Then, the determinant is a function of $a_1,~x_i,~y_i,~l,~m,$ and $n$; it can be written as: $$c_0+c_1x+c_2x^2+\cdots+c_6x^6,$$
where $c_i,~1\leq i\leq 6$, is a function on $a_1,~x_i,~y_i,~l,~m,$ and $n$. Next, the conditions in (6.1) are extended as follows, which leads to different signatures:
\begin{itemize}
	\item if $c_0<0$, sig=(5,1);
	\item if $c_0=0$ and $c_1<0$, sig=(4,1)
\end{itemize}

We apply the permutation groups $S_5$, $S_6$ and $S_6$ to the labels of the nodes in the Coxeter diagrams corresponding to $\mathcal{P}_3$, $\mathcal{P}_4$, and $\mathcal{D}_4$. This produces the desired data set consisting of all of the distinct Coxeter vectors as presented in Table \ref{table:sig344}. 

\begin{table}[h]
	{\footnotesize
		\begin{tabular}{c|c|c|c}
			\Xcline{1-4}{1.2pt}
			polytope &\quad \quad  $P^3$\quad  \quad   & \quad  \quad  $P^4$\quad  \quad   & \quad  \quad  $D^4$\quad  \quad  \\
			\hline
			\# SEILper matrix & 192 &26 &1200\\
			\hline
			$(3,1)$& 148 & 17   &  7 \\
			\hline
			$(4,1)$&148  & 17& 10\\
			\hline
		     total&148& 17 & 10 \\
		     \hline
		     \# distince Coxeter vector after permuting nodes &1122 &540& 864\\
			\Xcline{1-4}{1.2pt}
			
		\end{tabular}
	}
		
	\hspace*{0.5cm}
	\caption{Fushian and quasi-Fuchisian representations corresponding to $P^3$, $P^4$, and $D^4$. }
	\label{table:sig344}
\end{table}

In the following, we continue to use notations $\mathcal{P}_3$, $\mathcal{P}_4$ and $\mathcal{D}_4$ to denote the corresponding set of vectors. Firstly, we check all of the $109$ polytopes for the intersection types of $P^3$, $P^4$, and $D^4$ in their combinatorics. For example, regarding $P^3$, we are finding the five facets of $P_k$ such that the lattice formed by their intersection coincides with the lattice formed by the intersection of the five $2$-dimensional faces of the $3$-prism $P^3$. The polytopes having the intersection types of $P^3$, $P^4$, and $D^4$ in their combinatorics are listed in Table \ref{table:2potential}. Next, we use the data sets $\mathcal{P}^3$, $\mathcal{P}^4$, and $\mathcal{D}^4$ to filter the set of the SEILper matrices. Note that the three conditions for restriction types are ``saving" conditions, which foreshow their efficiency. The aforementioned ``saving" conditions,  $l4$- and $l5$-conditions, can be viewed as restrictions for the  intersection types of $3$-simplex and $4$-simplex, respectively.  The results after this procedure are arranged in column $6$ of Table \ref{table:2potential}.

\begin{table}[h]
	{\footnotesize
		\begin{tabular}{c|c|c|c|c|c|c}
			\Xcline{1-7}{1.2pt}
			\quad $d_k$ \quad &label of polytope &	\#potential & 	\#7	&intersection type (\#)	&\# (left)potential &	\quad  \#7 \quad \\
			\hline
			\multirow{2}{*}{5}&	254&	9467&	5&	$P^4$(5)&	42&	1\\
			\cline{2-7}
			    &311&	288	 &   2&	$P^4$(2)	&57	&1\\
			    \hline
			\multirow{4}{*}{4}&	249	&3910&	4&	$P^4$(2) &	40	&2\\
			\cline{2-7}
			    &255&	282	 &   3&	$P^4$(2) &	129	&2\\
			    \cline{2-7}
			   & 280&	148	&    3&	$P^4$(2)& 	0&	0\\ \cline{2-7}
			    &300&	383	&    2&	$P^3$(2) /$D^4$(1) &	336 / 6	&1\ \\ \hline
			\multirow{6}{*}{3}&	2&	160395&	6&	$P^3$(3)& 	139&	0\\ \cline{2-7}
			  & 10	&139143	&6&	$P^3$(3) / $D^4$(2) &	26915 / 27&	2\\ \cline{2-7}
		    	&242&	6744&	4&	$P^3$(2) / $P^4$(1) & 	324 / 40  &	1\\  \cline{2-7}
			   &247	&1546&	4&	$P^4$(1)&	545&	3\\ \cline{2-7}
		  	   &248&	2574&	4&	$P^4$(1)&	265&	3\\ \cline{2-7}
			  & 264&	514	  &  3&	$P^3$(2) / $D^4$(1) &	304 / 66&	3\\ \hline
			\multirow{7}{*}{2}&	4&	11420&	3&	$P^3$(4) &	397&	2\\ \cline{2-7}
			  &  9&	491	 &   2&	$P^3$(2) / $D^4$(1) &	486 / 204&	2\\ \cline{2-7}
			  & 18&	1058&	2&	$P^3$(2)&	309&	2\\ \cline{2-7}
			  & 61&	79	&    1&	$P^3$(2)&	44&	1\\ \cline{2-7}
			 &  68&	805	&    3&	$P^3$(2) / $D^4$(1) &	790 / 9	&2\\ \cline{2-7}
			  & 182&	2336&	4&	$P^3$(2)&	537&	3\\ \cline{2-7}
			  & 203&	908&	    4&	$P^3$(2) / $D^4$(1)& 	797 / 29&	3\\ 
			 \Xcline{1-7}{1.2pt}
			 	\multicolumn{7}{c}{}
			
		\end{tabular}
	}
	
		\hspace*{0.5cm}
		\hspace*{0.5cm}
		
	\caption{The ``intersection type (\#)" means there are \#-many combinatorial types of $P^3$, $P^4$, or $D^4$, possessed by the corresponding polytope. For example, ``$P^4(5)$" in the second line means there are $5$ sets of six facets of $P_{254}$ such that the lattice formed by them is a $4$-prism $P^4$.  The ``\# 7" means the maximal of the  numbers of $7$ of potential vectors in the corresponding set of SEILper potential vectors. It is obvious that the numbers of SEILper vectors decrease dramatically after this step. In addition, the values of \#7 decrease as well. }
	
	\label{table:2potential}
\end{table}

After deducing and modifying all of these limited SEILper potential matrices, we calculate (6.1) for each candidate. The sollowing logic is adopted:
\begin{enumerate}
	\item Prepare conditions (6.1) for a given potential vector
	\item Calculate the Gr\"{o}bner bases of the polynomials involved and modify the condition (6.1) correspondingly.
	\item Solve the equations and inequality in (2) simultaneously and restrict the running time to a maximum of $30s$.
	\item If the calculation does not finish in $30s$, the program aborts and turns to solve the equations and inequality in (1) directly.
\end{enumerate}

Although we set a threshold of $30s$, for almost all of the vectors, the calculations can be done within $1s$. Moreover, we use the logic above (1)-(4) rather than continuing to wait for the solutions of the equations and inequality in (2). This is because the function  \emph{GroebnerBasis} in \emph{Mathematica} either works quite quickly or consumes an unaffordable amount of time. Alternatives are therefore needed. The last step is to check whether the signature is indeed $(5,1)$ and whether $\displaystyle \frac{\pi}{\arccos~ (y_i/2)}$ is an integer. 

After conducting all of these procedures in \emph{Mathematica}, we find that only six of all simple $5$-polytopes with $9$ facets admit compact hyperbolic structure. We also glue the $4$-prism to the one with orthogonal 4-prism ends. The results are reported in Table \ref{table:result}. The polytopes with labels in red are the ``basis" polytopes, and the ones on last line can be obtained by gluing prism ends to those of the second line from the bottom. The Coxeter diagrams and length information are shown in the end (pages 23--26). 
 
\begin{table}[h]
	{\footnotesize
		\begin{tabular}{c|c|c|c|c|c|c}
			\Xcline{1-7}{1.2pt}
			$d_k$ &	\multicolumn{3}{c|}{6} & $5$&	$4$ &	$3$\\
			\hline
			polytope&\quad	{\color{red}312}\quad&	\quad {\color{red}319}\quad	& \quad {\color{red}322}\quad&	\quad {\color{red}302}\quad&\quad	{\color{red}313}\quad&	\quad 284\quad\\
			\hline
			\# (selected) SEILper potential	&1	&3&	5&	5	&68&	60\\
			\hline
			\# gram (of basis vectors) &	1&	3&	5&	1&	1&	1\\
			\hline
			\# gram after suitably gluing & \multirow{2}{*}{1} &\multirow{2}{*}{22}& \multirow{2}{*}{18}&\multirow{2}{*}{6}&\multirow{2}{*}{3} &\multirow{2}{*}{N}\\
			
			$4$-prisms to the orthoganal ends&&&&&&\\	
			
			\Xcline{1-7}{1.2pt}
			
		\end{tabular}
	}
	
	\vspace{0.5cm}
	\caption{Results of the compact hyperbolic Coxeter $5$-polytopes with 9 facets. The value of polytope labeled by $284$ on the last line is ``N", which means the $l_5\_basis$ set of $P_{284}$ is empty and we are not be able to glue it with any prism end.}
	\label{table:result}
\end{table}

\newpage
\section{Validation and Results}\label{section:vadilation}
1. ``Basis approach" vs. ``direct approach".

We calculate SEILper potential matrices without using  $l_5\_basis$-conditions for those polytopes having $4$-simplex facets. The restriction sets regarding the intersection types of $P^3$, $P^4$, and $D^4$ are then applied. After checking for the signature obstructions, we  finally obtain the Gram matrices corresponding to all of the possible compact hyperbolic polytopes. The results are reported in Table \ref{table:direct}. They are the same as the previous work generated via ``basis approach".

\begin{table}[H]
	{\footnotesize
		\begin{tabular}{c|ccccccccccc}
			\Xcline{1-12}{1.2pt}
			
			grp	&\multicolumn{3}{c|}{6}		&\multicolumn{6}{c}{5}			&&		\\
			\Xcline{1-12}{1.2pt}
			polytopes&	312&	319	&\multicolumn{1}{c|}{322}&	&302&	311&	315&	318&	320&	321&\\
			\#SEILper&	1&	22&	\multicolumn{1}{c|}{18}&&	18&	560&	139	&2	&12	&9&\\
			$P^3$/$P^4$/$D^4$ restrictions&	N&	N&	\multicolumn{1}{c|}{N}&&	N&	99&	N&	N&	N&	N&\\
			\#Gram&	1&	22&	\multicolumn{1}{c|}{18}&&	6&	0&	0&	0&	0&	0&\\
			\Xcline{1-12}{1.2pt}
			grp&	\multicolumn{9}{c}{4}&&		\\	
			\Xcline{1-12}{1.2pt}					
			polytopes&&	295&	300&	301&	303&	310&	313	&314&	316	&317&\\
			\#SEILper&&	0&	709&	15&	0&	0&	110&	0&	14&	0&\\
			P3/P4/D4 restrictions&&	N&	586/8&	N&	N&	N&	N&	N&	N&	N&\\
			\#Gram&&	0&	0&	0&	0&	0&	3&	0&	0&	0&\\
	
			\Xcline{1-12}{1.2pt}
			grp&	\multicolumn{11}{c}{3}\\	
			\Xcline{1-12}{1.2pt}							
			polytopes&	294&	296	&297&	298&	299&	304	&305&	306&	307&	308&	309\\
			\#SEILper&	83	&37	&0	&11&	0&	0&	0&	0&	0&	96&	0\\
			P3/P4/D4 restrictions&	N&	N&	N&	N&	N&	N&	N&	N&	N&	N&	N\\
			\#Gram&	0&	0&	0&	0	&0&	0&	0&	0&	0&	0&	0\\
			\Xcline{1-12}{1.2pt}

		\end{tabular}
	}
	
	\hspace*{0.5cm}
	\caption{Results obtained by direct approach. }
	\label{table:direct}
\end{table}

2. A. Burcroff presents $50$ polytopes in \cite{Amanda:2022}. A careful check reveals that we have covered all of these polytopes and obtained one more. The correspondence of notions of the polytopes, which admit a hyperbolic structure, between our list and Burcroff's are presented in Table \ref{comparison}.

\begin{table}[h]
	{\footnotesize
		\begin{tabular}{c|cccccc}
			\Xcline{1-7}{1.2pt}

MZ&	322	&319&	312&	302&	313&	284\\
\hline
A. Burcroff&	$H_2$&	$H_1$	&&	$H_3$&	$H_4$&	$H_5$\\
\Xcline{1-7}{1.2pt}

\end{tabular}
}

\hspace*{0.5cm}
\caption{Notion correspondence between our list and the list in \cite{Amanda:2022}. }
\label{comparison}
\end{table} 

\newpage
\newgeometry{left=0.5cm,right=0.5cm,top=1.5cm,bottom=2cm}

1. Coxeter diagrams for $P_{322}$.

\vspace{0.5cm}

\begin{figure}[H]
	\scalebox{0.3}[0.3]{\includegraphics {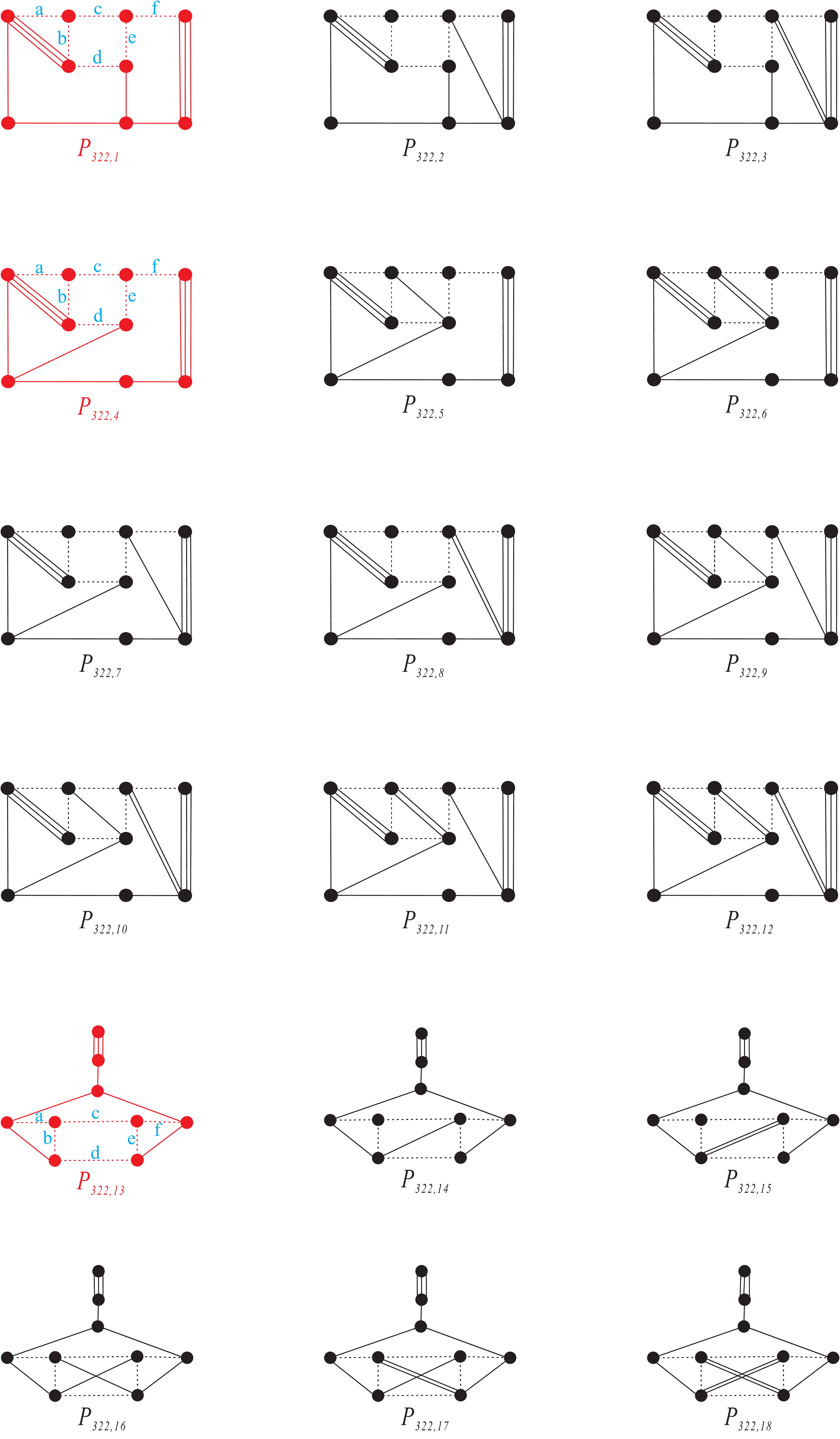}}
	\caption{The $18$ compact hyperbolic Coxeter $5$-polytopes with $9$ facets over polytope   $P_{322}$.} \label{figure:p322}
\end{figure}

\newpage
2. Coxeter diagrams for $P_{319}$ (part 1).

\vspace{0.5cm}

\begin{figure}[H]
	\scalebox{0.3}[0.3]{\includegraphics {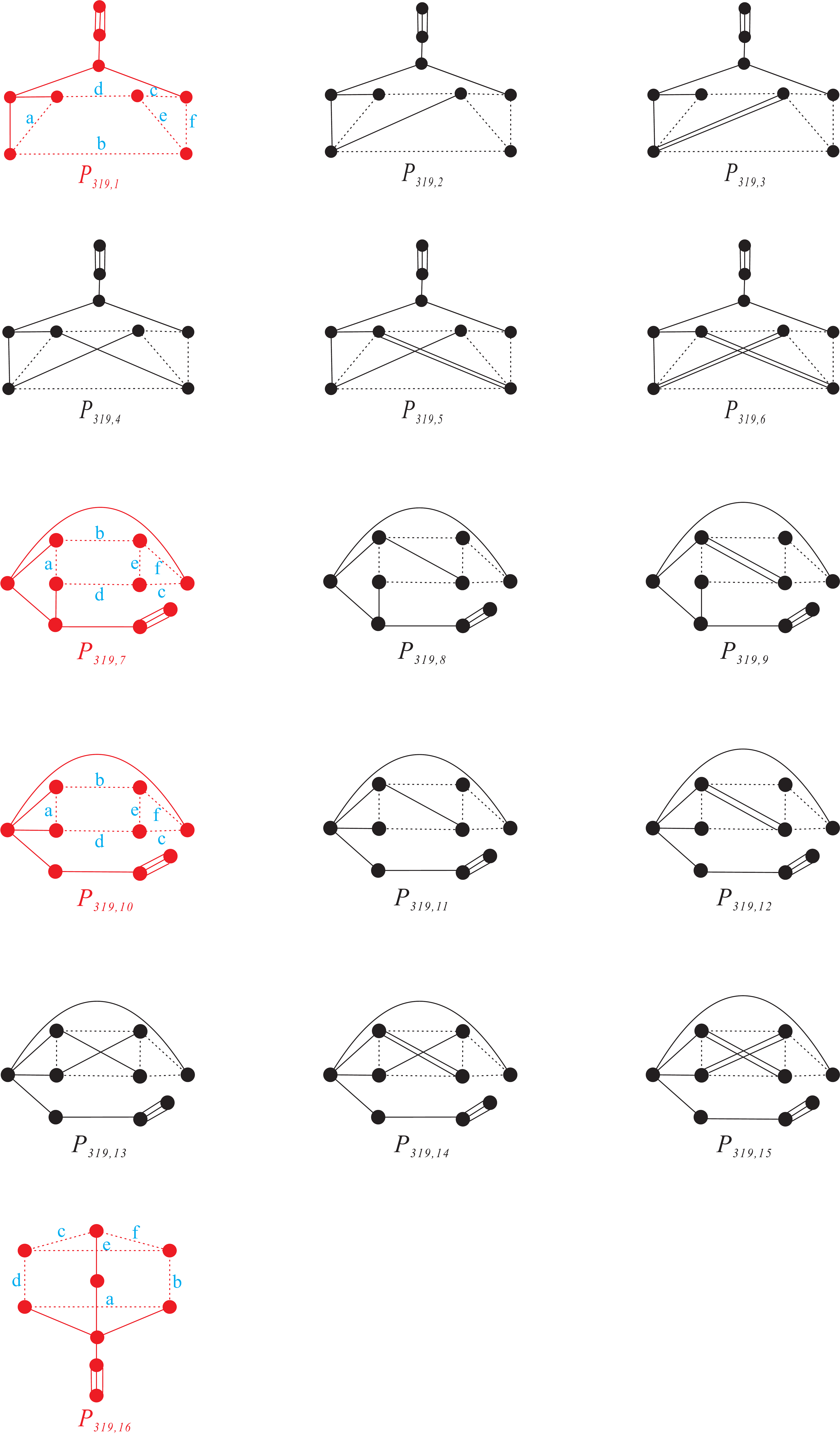}}
	\caption{$16$ of the $22$ compact hyperbolic Coxeter $5$-polytopes with $9$ facets over polytope  $P_{319}$.} \label{figure:p319}
\end{figure}

\newpage
3. Coxeter diagrams for $P_{319}$ (part 2), $P_{312}$, $P_{302}$, $P_{313}$.

\vspace{0.5cm}

\begin{figure}[H]
	\scalebox{0.3}[0.3]{\includegraphics {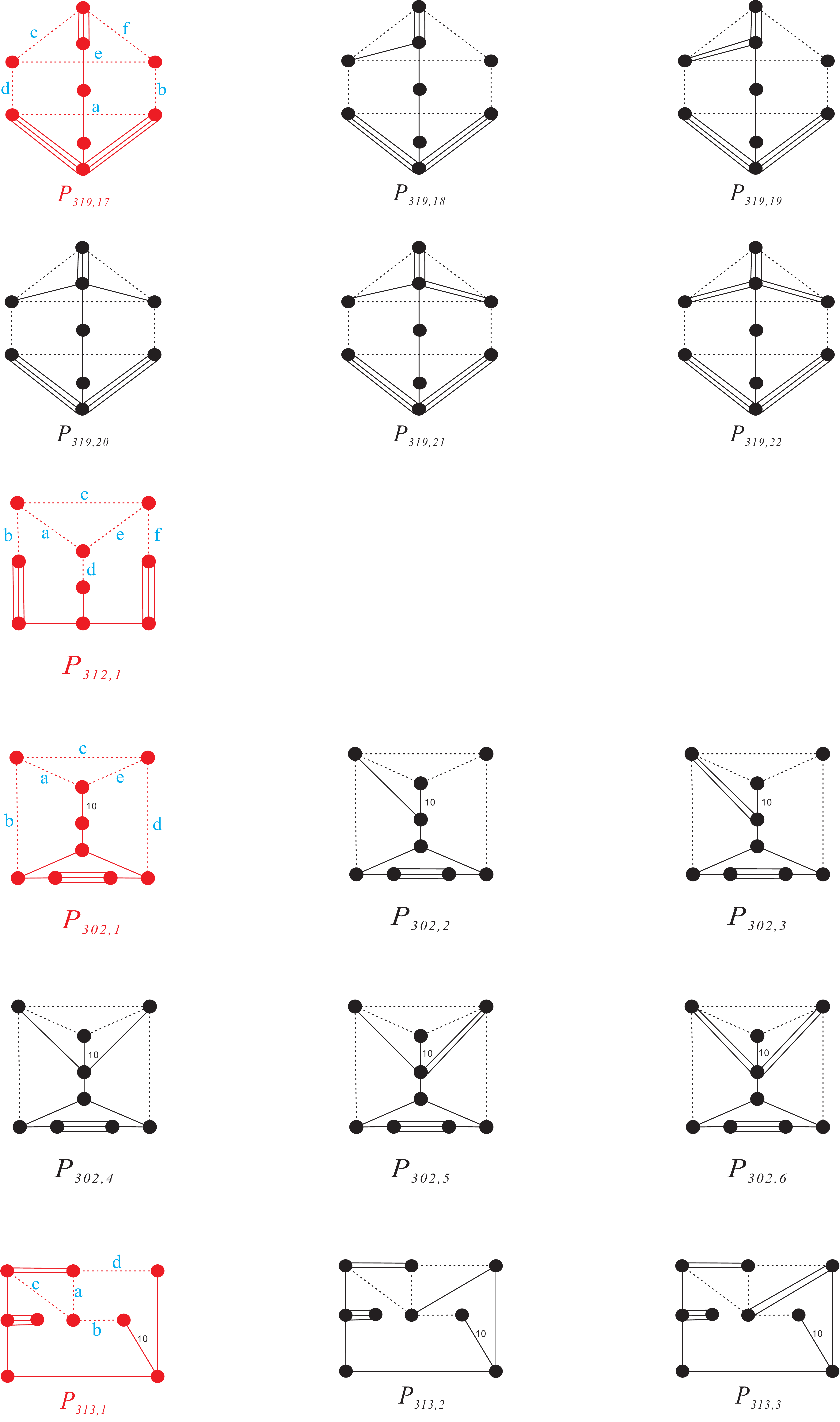}}
	\caption{$6$ of $22$, $1$, $6$, $3$ compact hyperbolic Coxeter $5$-polytopes with $9$ facets over polytopes  $P_{319}$, $P_{312}$, $P_{302}$, $P_{313}$, respectively.} \label{figure:p19120213}
\end{figure}

\newpage
\newgeometry{left=0.5cm,right=0.5cm,top=1cm,bottom=1.5cm}
4. Coxeter diagram for $P_{284}$.

\begin{figure}[H]
	\scalebox{0.3}[0.3]{\includegraphics {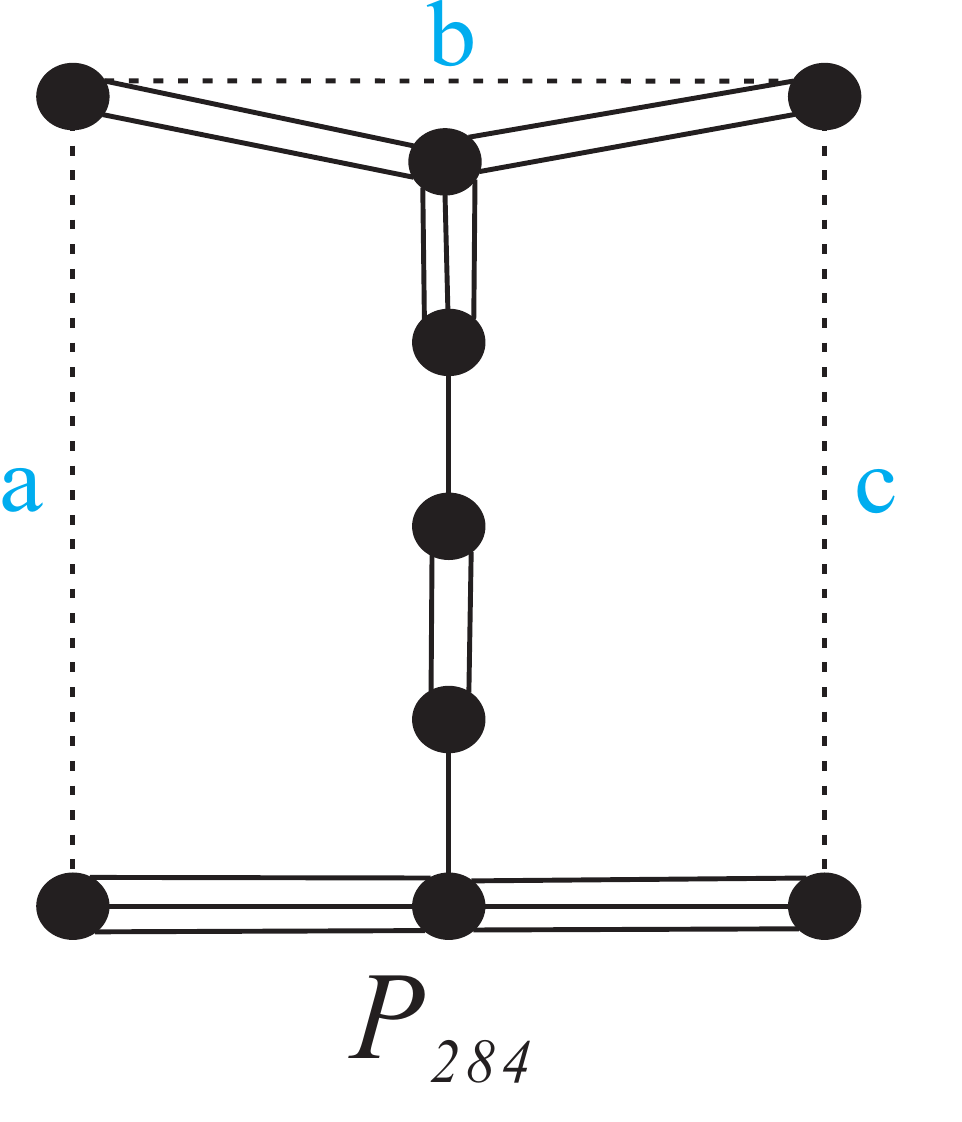}}
	\caption{The only compact hyperbolic Coxeter $5$-polytopes with $9$ facets over polytope   $P_{284}$} \label{figure:p284}
\end{figure}

	\begin{table}[H]
		\renewcommand{\arraystretch}{1.9}
		\resizebox{17cm}{!}{ %
			\begin{tabular}{|l|l|l|l|}
				\hline
				\multirow{2}{*}{} &a  &b  &c  \\ \cline{2-4} 
				&d  &e  &f  \\ \hline
				\multirow{2}{*}{$P_{322,1}$} &$\frac{1}{2}\sqrt{\frac{1}{2}(6+\sqrt{5})}$ & $\frac{1}{4}\sqrt{23+8\sqrt{5}+\sqrt{5(63+26\sqrt{5})}}$ & $\frac{1}{4}\sqrt{946+423\sqrt{5}+\sqrt{1749835+782550\sqrt{5}}}$  \\ \cline{2-4} 
				&$\frac{1}{4}(5+2\sqrt{5}+\sqrt{63+26\sqrt{5}})$ & $\frac{1}{2}\sqrt{\frac{1}{2}(1257+562\sqrt{5}+3\sqrt{349967+156510\sqrt{5}}}$  &$\frac{1}{2}\sqrt{\frac{1}{2}(9+\sqrt{5})}$  \\ \hline
				\multirow{2}{*}{$P_{322,4}$} &  $\frac{1}{2}(1+\sqrt{5})$ & $\frac{1}{4}(1+\sqrt{5}+2\sqrt{8+3\sqrt{5}})$ & $\frac{9}{2}+2\sqrt{5}+\frac{1}{4}\sqrt{1382+618\sqrt{5}}$  \\ \cline{2-4} 
				& $\frac{1}{4}(3+\sqrt{5}+2\sqrt{8+3\sqrt{5}})$ &  $\frac{1}{2}(9+4\sqrt{5}+\sqrt{132+59\sqrt{5}}) $ & $\frac{1}{2}\sqrt{\frac{1}{2}(9+\sqrt{5})}$ \\ \hline
				\multirow{2}{*}{$P_{322,13}$} & $\sqrt{\frac{23}{8}+\frac{9\sqrt{5}}{8}}$ &   $\frac{1}{2}\sqrt{6+\sqrt{5}}$  &  $\frac{1}{4}(15+7\sqrt{5})$  \\ \cline{2-4} 
				&$\frac{1}{4}(3+\sqrt{5})$  &  $\frac{1}{2}\sqrt{6+\sqrt{5}}$  & $\sqrt{\frac{23}{8}+\frac{9\sqrt{5}}{8}}$ \\ \hline
				\multirow{2}{*}{$P_{319,1}$} &  $\frac{1}{4}(2+\sqrt{5})$ &   $\frac{1}{4}\sqrt{\frac{5}{2}(23+9\sqrt{5})} $  & $\sqrt{\frac{23}{8}+\frac{9\sqrt{5}}{8}}$  \\ \cline{2-4} 
				& $\frac{1}{4}\sqrt{\frac{5}{2}(23+9\sqrt{5})} $  & $\frac{1}{8}(7+5\sqrt{5})$  &$\sqrt{\frac{23}{8}+\frac{9\sqrt{5}}{8}}$  \\ \hline
				\multirow{2}{*}{$P_{319,7}$} &  $\frac{1}{4}(5+3\sqrt{5}+2\sqrt{\frac{57}{2}+\frac{25\sqrt{5}}{2}})$ &   $\frac{5+3\sqrt{5}+\sqrt{46+18\sqrt{5}}}{4\sqrt{2}} $  &  $\frac{1}{2}(1+\sqrt{5})$  \\ \cline{2-4} 
				& $\frac{5}{2}+\sqrt{5}+\frac{1}{4}\sqrt{114+50\sqrt{5}}$  & $\frac{1}{4}\sqrt{47+17\sqrt{5}+2\sqrt{570+250\sqrt{5}}}$ &$\frac{1}{2}\sqrt{\frac{1}{2}(6+\sqrt{5})}$  \\ \hline
				\multirow{2}{*}{$P_{319,10}$} &  $2+\sqrt{5}$ &  $3+\sqrt{5} $ &  $\frac{1}{2}(1+\sqrt{5})$  \\ \cline{2-4} 
				& $3+\sqrt{5}$ & $2+\sqrt{5}$  &$\frac{1}{2}(1+\sqrt{5})$  \\ \hline
				\multirow{2}{*}{$P_{319,16}$} & $\frac{1}{4}(19+9\sqrt{5})$ &   $\frac{1}{4}\sqrt{5(57+25\sqrt{5})}$ &  $\frac{1}{2}\sqrt{\frac{1}{2}(6+\sqrt{5})}$  \\ \cline{2-4} 
				& $\frac{1}{4}\sqrt{5(57+25\sqrt{5})}$  &  $\frac{1}{8}(7+5\sqrt{5})$ & $\frac{1}{2}\sqrt{\frac{1}{2}(6+\sqrt{5})}$ \\ \hline
				\multirow{2}{*}{$P_{319,17}$} & $\frac{1}{4}(5+\sqrt{5})$ &   $\sqrt{\frac{691}{8}+\frac{309\sqrt{5}}{8}}$ &  $\frac{1}{2}\sqrt{\frac{1}{2}(9+\sqrt{5})}$  \\ \cline{2-4} 
				& $\sqrt{\frac{691}{8}+\frac{309\sqrt{5}}{8}}$  & $\frac{1}{4}(119+55\sqrt{5})$  &$\frac{1}{2}\sqrt{\frac{1}{2}(9+\sqrt{5})}$  \\ \hline
				\multirow{2}{*}{$P_{312,1}$} & $\frac{1}{2}\sqrt{\frac{5}{2}+\sqrt{5}}$ &  $\frac{1}{4}\sqrt{15+\sqrt{5}}$ & $1+\frac{\sqrt{5}}{2}$   \\ \cline{2-4} 
				& $\frac{1}{2}\sqrt{4+\sqrt{5}}$  & $\frac{1}{2}\sqrt{\frac{5}{2}+\sqrt{5}}$  &$\frac{1}{4}\sqrt{15+\sqrt{5}}$   \\ \hline
				
				\multirow{2}{*}{$P_{302,1}$} &$\frac{1}{2}\sqrt{4+\sqrt{5}}$  & $\frac{1}{2}\sqrt{30+13\sqrt{5}}$  & $\frac{1}{4}(31+15\sqrt{5})$  \\ \cline{2-4} 
				& $\frac{1}{2}\sqrt{30+13\sqrt{5}}$   & $\frac{1}{2}\sqrt{4+\sqrt{5}}$  &  \\ \hline
				
				\multirow{2}{*}{$P_{313,1}$} &$\frac{1}{2}\sqrt{3+\sqrt{5}}$  & $\frac{1}{2}\sqrt{5+\sqrt{5}}$  & $\frac{1}{2}(3+\sqrt{5})$  \\ \cline{2-4} 
				& $\sqrt{5+2\sqrt{5}}$   &   &  \\ \hline
				
				$P_{284}$ &  $\frac{1}{2}(1+\sqrt{5})$ &  $\frac{1}{2}(1+\sqrt{5})$ & $\frac{1}{2}(1+\sqrt{5})$ \\ \hline
			\end{tabular}
		}
	
		\caption{The $\cosh \rho_{ij}$, denoted by letter alphabetically, of the result polytopes  .}
		\label{table: length}
	\end{table}

\newpage
\restoregeometry

\end{document}